\documentclass{article}
\usepackage{amssymb}
\usepackage{amsmath}

\newtheorem{thm}{Theorem}[section]
\newtheorem{cor}[thm]{Corollary}
\newtheorem{lem}[thm]{Lemma}
\newtheorem{prop}[thm]{Proposition}
\newtheorem{defn}[thm]{Definition}
\newtheorem{rem}[thm]{Remark}

\newcommand{\norm}[1]{\left\Vert#1\right\Vert}
\newcommand{\Real}{\mathbb R}

\newcommand{\half}{\frac{1}{2}}
\newcommand{\diffd}{\mathrm{d}} 
\newcommand{\Prob}{\mathbb{P}} 
\newcommand{\Expct}{\mathbb{E}} 

\newcommand{\Proof}{\noindent\textbf{Proof: }}
\newcommand{\Remark}{\textbf{Remark: }}

\newcommand{\F}{\mathcal{F}} 
\newcommand{\eL}{\mathcal{L}} 
\newcommand{\N}{\mathbb{N}} 
\newcommand{\Ind}{\mathbf{1}} 

\newcommand{\cond}[2]{\bigl(#1\lvert#2\bigr)}
\newcommand{\ccond}[2]{\Bigl(#1\bigl\vert#2\Bigr)}
\newcommand{\seq}[1]{\bigl(#1\bigr)}  
\newcommand{\bseq}[1]{\Bigl(#1\Bigr)}  
\newcommand{\bbseq}[1]{\biggl(#1\biggr)}  
\newcommand{\bsseq}[1]{\Bigl[#1\Bigr]}  
\newcommand{\set}[1]{\bigl\{#1\bigr\}}
\newcommand{\bset}[1]{\Bigl\{#1\Bigr\}}
\newcommand{\bbset}[1]{\biggl\{#1\biggr\}}
\newcommand{\abs}[1]{\bigl\vert#1\bigr\vert}

\newcommand{\ip}[2]{\bigl\langle #1,#2 \bigr\rangle} 

\newcommand{\node}[1]{\text{node}_{#1}}
\newcommand{\Q}{\mathbb{Q}} 

\newcommand{\insidewith}{(\tilde{\lambda}(\theta),0]}

\newcommand{\outsidewithout}{<\tilde{\lambda}(\theta)}
\newcommand{\outsidewith}{\leq\tilde{\lambda}(\theta)}

\newcommand{\ie}{i.e.\ }

\newcommand{\qed}{\hfill$\Box$\bigskip}

\ifx\pdfoutput\undefined 
\else
\RequirePackage[colorlinks,hyperindex,plainpages=false]{hyperref}
\fi

\begin{document}

\title{Spine proofs for $\eL^p$-convergence of branching-diffusion \\ martingales\thanks{
This \texttt{arXiv} article is an revision of \emph{Spine proofs for
$\eL^p$-convergence of branching-diffusion martingales}, (2004),
no.~0405, Mathematics Preprint, University of Bath. }}
\author{Robert Hardy and Simon C. Harris\thanks{
Email: \texttt{S.C.Harris@bath.ac.uk},
Web: http://people.bath.ac.uk/massch}\\ University of Bath, UK}
\date{\today}
\maketitle

\begin{abstract}
Using the foundations laid down in Hardy and Harris \cite{spine_foundations},
we present new spine proofs of the $\eL^p$-convergence $(p\geq1)$ of some key
`additive' martingales for three distinct models of branching diffusions,
including new results for a multi-type branching Brownian motion
and discussion of left-most particle speeds.
The spine techniques we develop give clear and simple arguments in the spirit of
the conceptual spine proofs found in Kyprianou \cite{kyprianou:travelling_waves} and Lyons \emph{et al}
\cite{lyons_pemantle_peres:conceptual_LlogL,lyons:simple_biggins_convergence,
lyons_kurtz_pemantle_peres:conceptual_Kesten_Stigum}, and they
should also extend to more general classes of branching diffusions.
Importantly, the techniques in this paper also pave the way for the path large-deviation results for
branching diffusions found in Hardy and Harris \cite{bbm_large_deviations,preprint-HW_large_deviations}.
\end{abstract}

\section{Overview}
In this article we use a change of measure together with spine
techniques to analyze the $\eL^p$-convergence properties (for $p\geq 1$)
of the strictly-positive `additive' martingales for three different
models of branching diffusions.
It is a common feature of these diffusion models,
where there is actually a family of such martingales $\set{Z_\lambda: \lambda \in \Real}$,
that for all $\lambda$ within an open interval about the origin the martingale
$Z_\lambda$ is convergent in $\eL^p$ for some $p\geq 1$
subject to a suitable $p^{\textrm{th}}-$moment ($p>1$) or $L\log L$ ($p=1$) condition on the offspring distribution;
for $\lambda$ outside of this interval, or if the $L\log L$ condition fails for the offspring distribution,
the limit of $Z_\lambda$ is almost surely null.

The first model we consider is a branching Brownian motion (BBM) with random family sizes.
After introducing the fundamental `additive' $Z_\lambda$ martingales
and describing a `spine' construction for the BBM under a change of measure using $Z_\lambda$,
we will recall Kyprianou's \cite{kyprianou:travelling_waves}  $\eL^1$-convergence result
before stating necessary and sufficient conditions for $\eL^p$-convergence of the $Z_\lambda$ martingales.
Our new proof of martingale $\eL^p$-convergence for BBM uses `spine' change of measure techniques
and, early on in Section \ref{chapter_neveu_lp_alternative_proof},
we will include a summary of the underlying
space and filtrations that we shall use for our spine
techniques throughout this paper; a foundation
article \cite{spine_foundations} contains full details,
but we have tried to keep this article reasonably self-contained.
Note that for BBM,
the spine construction first appeared in Chauvin and Rouault \cite{chauvin:KPP_app_to_spatial_trees},
whilst Kyprianou \cite{kyprianou:travelling_waves} really exploited `spine' methods in his proofs,
however our spine approach does possess some significant differences from others.
Neveu\cite{neveu} used classical techniques
for $\eL^p$-convergence in the special case of binary branching.
Also see Harris \cite{art:twFKPP} for further discussion of martingale convergence in BBM and applications.

In Section \ref{chapter_MTBBM}, we look at a finite-type BBM model where the type of each particle
controls the rate of fissions, the offspring distribution and the spatial diffusion.
First, we will extend Kyprianou's \cite{kyprianou:travelling_waves} approach to give the
analogous $\eL^1$-convergence result for this multi-type BBM model.
We will also briefly discuss the rate of convergence of the martingales to zero and the
speed of the spatially left-most particle within the process.
Next, we give a new result on $\eL^p$-convergence criteria, extending our earlier spine
based proof developed for the single-type BBM case.

The third model of Section \ref{chapter:HW_model}
has a \emph{continuous}-type-space where the type
of each particle moves independently as an Orstein-Uhlenbeck
process on $\Real$. This branching diffusion was first introduced
in Harris and Williams \cite{art:largedevs1} and has also been
investigated in Harris \cite{harris:gibbs_boltzmann}, 
Git \emph{et al} \cite{ghh}
and Kyprianou and Engl\"{a}nder
\cite{kyprianou_englander}.

Proofs for each of these models each run along similar lines and the techniques are
quite general, and it is a powerful feature of the spine approach that this is possible.
More classical techniques based on the expectation semigroup are simply
not able to generalize easily, since they often require either
some \emph{a priori} bounds on the semigroup or involve difficult
estimates -- for example, in Harris and Williams
\cite{art:largedevs1} their important bound of a non-linear term
is made possible only by the existence of a good $\eL^2$ theory
for their operator, and this is not generally available.

Briefly, to prove that the martingale converges in $\eL^p$ for
some $p > 1$ we use Doob's theorem, and therefore need only to
show that the martingale is \emph{bounded} in $\eL^p$. The
\emph{spine decomposition} is an excellent tool here for showing
boundedness of the martingale since it reduces difficult
calculations over the whole collection of branching particles to
just the single spine process.
We find the same conditions are also necessary for $\eL^p$-boundedness of the martingale
when $p>1$ by just considering the contributions along the spine at times of fission
and observing when these are unbounded.
Otherwise, to determine whether the martingale is merely $\eL^1$-convergence or
has an almost-surely zero limit,
we determine whether the martingale is almost-surely bounded or not
under its own change of measure
-- this was Kyprianou's \cite{kyprianou:travelling_waves}
approach and relies on a measure-theoretic
result and has become standard in the spine
methodology since the important work of Lyons \emph{et al}
\cite{lyons_pemantle_peres:conceptual_LlogL,
lyons:simple_biggins_convergence,
lyons_kurtz_pemantle_peres:conceptual_Kesten_Stigum}.
Spine and size-biasing techniques have already proved extremely useful in many other branching process
situations, for example,
also see Athreya \cite{athreya:bernoulli_change_of_measure_MC},
Biggins and Kyprianou \cite{bigginkypMTBP},
Geiger \cite{geiger:splitting_trees, Geiger:Sizebiasedsplittingtrees,
 geiger:elementary_GW_proofs, Geiger:PoissonGWtrees,  Geiger:infvartree},
Georgii and Baake \cite{georgii},
Iksanov \cite{iksanov},
Olofsson \cite{olofsson:xlogx}, Rouault and Liu \cite{rouault_liu:two_measures_on_tree_boundary}
and Waymire and Williams \cite{Waymire},
to name just a few.

There are a number of reasons why we may be interested in knowing
about the $\eL^p$ convergence of a martingale: in Neveu's original
article \cite{neveu} it was a means to proving
$\eL^1$-convergence of martingales which can then be used to represent (non-trivial) travelling-wave
solutions to the FKPP reaction-diffusion equation as well as in understanding the growth and spread of the BBM,
whilst
Git \emph{et al} \cite{ghh} and Asmussen and Hering
\cite{asmussen_hering:strong_limits_no_immigration} have used it
to deduce the almost-sure rate of convergence of the martingale to
its limit.
Of equal importance are the \emph{techniques} that we
use here:
the convergence of other additive martingales can be
determined with similar techniques, for example, see an application to
a BBM with quadratic breeding potential in J.W.Harris and S.C.Harris \cite{hhx2};
similar ideas have also been used in proving a lower
bound for a number of problems in the large-deviations theory of
branching diffusions -- we have used the spine decomposition with
Doob's submartingale inequality to get an upper-bound for the
growth of the martingale under the new measure which then leads to
a lower-bound on the probability that one of the diffusing
particles follows an unexpected path -- see Hardy and Harris
\cite{bbm_large_deviations} for a spine-based proof of
the large deviations principle for branching Brownian motion, and
see Hardy and Harris \cite{preprint-HW_large_deviations} for a proof of a lower
bound in the model that we consider in section \ref{chapter:HW_model}.

\section{Branching Brownian motion}\label{chapter_neveu_lp_alternative_proof}

Consider a branching Brownian motion (BBM) with constant branching
rate $r$, which is the branching process whereby particles diffuse
independently according to a (driftless) Brownian motion and at any moment
undergo fission at a rate $r$ to produce a random number of offspring,
$1 + A$, where $A$ is an independent random variable with distribution
\[
P(A = i) = p_i,\qquad i \in \set{0,1,\ldots},
\]
such that $m := P(A) = \sum_{i=0}^\infty i \, p_i < \infty$.
Offspring move off from their parents point of fission, and
continue to evolve independently as above.
We suppose that the probabilities of this process are $\set{P^x: x \in \Real}$ so
that $P^x$ is a measure defined on the natural filtration
$(\F_t)_{t \geq 0}$ such that it is the law of the process
initiated from a single particle positioned at $x$.

Suppose that the configuration of this branching Brownian motion
at time $t$ is given by the $\Real$-valued point process
$\mathbb{X}_t := \set{X_u(t): u \in N_t}$ where $N_t$ is the set
of individuals alive at time $t$. It is well known that
for any $\lambda \in \Real$,
\begin{equation}
Z_\lambda(t) := \sum_{u \in N_t} e^{-rmt}e^{\lambda X_u(t) - \half
\lambda^2 t} = \sum_{u \in N_t} e^{\lambda X_u(t) - E_\lambda t}
\end{equation}
where $E_\lambda := -\lambda c_\lambda:= \half\lambda^2 + rm$
defines a \emph{positive} $P$-martingale, so
$Z_\lambda(\infty) := \lim_{t \to \infty} Z_\lambda(t)$ is almost
surely finite under each $P^x$.

We are going to use a change of measure together with the
so-called spine decomposition to determine the conditions under
which this martingale is $\eL^p(P)$-convergent for some $p> 1$.

\begin{thm}\label{thm:BBM_c-o-m_with_Z}
If we define the measure $\Q_\lambda^x$ via
\begin{equation}
\frac{\diffd \Q_\lambda^x}{\diffd P^x}\bigg\vert_{\F_t} =
\frac{Z_\lambda(t)}{Z_\lambda(0)} = e^{-\lambda x} Z_\lambda(t),
\label{BBM_Z_c-o-m}
\end{equation}
then it follows that under $\Q_\lambda^x$ the point process
$\mathbb{X}_t$ evolves as follows:
\begin{itemize}
\item starting from position $x$, the original ancestor diffuses
according to a Brownian motion on $\Real$ with drift $\lambda$;

\item at an accelerated rate $(1+m)r$ the particle undergoes
fission producing $1+\tilde{A}$ particles, where the distribution
of $\tilde{A}$ is independent of the past motion but is
\emph{size-biased}:
\[
\tilde{\Q}_\lambda(\tilde{A} = i) = \frac{(i+1)p_i}{m+1},\qquad i
\in \set{0,1,\ldots}.
\]

\item with equal probability, one of these offspring particles is
selected;

\item this chosen particle repeats stochastically the behaviour of
the parent with the size-biased offspring distribution;

\item the other particles initiate, from their birth position, an
independent copy of a $P^\cdot$ branching Brownian motion with
branching rate $r$ and family-size distribution given by $A$
(which is without the size-biasing).
\end{itemize}
\end{thm}
In this construction, the individuals that are selected to have a
drift of $\lambda$ make up a (random) line of descent which has
come to be referred to as the \emph{spine}.
The phenomena of size-biasing along the spine is a common feature of such measure changes
when random offspring distributions are present.
This \emph{spine} construction for BBM can also be seen in
Chauvin and Rouault \cite{chauvin:KPP_app_to_spatial_trees},
Kyprianou \cite{kyprianou:travelling_waves} and
J.Harris, S.C.Harris and Kyprianou \cite{hhk}.

In particular, Kyprianou \cite{kyprianou:travelling_waves}
used this change of measure and other spine
techniques to give necessary and sufficient conditions for
$\eL^1$-convergence of the $ Z_\lambda$ martingales.
By a natural symmetry, without loss of generality we will throughout suppose that $\lambda\leq 0$, then:
\begin{thm}\label{BBM_cgce_theorem}
Let $\tilde{\lambda} := -\sqrt{2rm}$ so that $c_\lambda:=-E_\lambda/\lambda$ attains local maximum
at $\tilde{\lambda}$.
For each $x\in \Real$, the limit $Z_\lambda(\infty) := \lim_{t \to \infty} Z_\lambda(t)$ exists
$P^x$-almost surely where:
\begin{itemize}
\item  if $\lambda \leq \tilde{\lambda}$ then $Z_\lambda(\infty) = 0$ $P^x$-almost surely;
\label{kyp_part_1}
\item  if $\lambda \in (\tilde{\lambda},0]$ and $P(A \log^+A)=\infty$ then $Z_\lambda(\infty) = 0$ $P^x$-almost surely;
\label{kyp_part_2}
\item if $\lambda \in (\tilde{\lambda},0]$ and $P(A \log^+A)<\infty$ then $Z_\lambda(t)\rightarrow Z_\lambda(\infty)$
almost surely and in $\eL^1(P^x)$.
\label{kyp_part 3}
\end{itemize}
\end{thm}

We will give a generalisation of this result to a multi-type BBM in the next section,
proving it there by extending Kyprianou's proof for BBM as well as discussing the rate of convergence
to zero of the martingales and left-most particle speeds.
In fact, in many cases where the martingale has a non-trivial limit,
the convergence will also be much stronger than merely in $\eL^1(P^x)$, as
indicated by the following $\eL^p$-convergence result:
\begin{thm}\label{BBM_L_p_theorem}
For each $x \in \Real$, and for each $p \in (1, 2]$:
\begin{itemize}
\item $Z_\lambda(t)\rightarrow Z_\lambda(\infty)$ almost surely and in $\eL^p(P^x)$ if
$p\lambda^2 < 2mr$ and $P(A^p)<\infty$
\label{Lp_part_1}
\item $Z_\lambda$ is unbounded in $\eL^p(P^x)$, that is, $P^x(Z_\lambda(\infty)^p)=\infty$ if
$p\lambda^2 > 2mr$ or $P(A^p)=\infty$.
\label{Lp_part_2}
\end{itemize}
\end{thm}

We shall give a spine-based proof of this $\eL^p$-convergence theorem
in Section \ref{sec:proof_of_BBM}, but also
see Neveu \cite{neveu} for sufficient conditions in the special case of binary branching
at unit rate using more classical techniques.
Iksanov \cite{iksanov} also uses similar spine techniques in the study of the
branching random walk.

\subsection{The underlying space and filtrations}\label{sec_BBM:foundations}

For a clear understanding of the spine techniques that we shall
use in all our models, we need a more precise description of
spines than the pathwise construction given in Theorem
\ref{thm:BBM_c-o-m_with_Z}. We give fuller details in Hardy and Harris
\cite{spine_foundations}, but to make this article
self-contained we now briefly lay out the principal elements. The
reader who is familiar with the work of Lyons \emph{et al} or
Kyprianou \cite{kyprianou:travelling_waves} will
notice significant differences in our approach via our use of the
filtrations on the single underlying space.

The basic ideas of our approach is quite straightforward:
given the original BBM, we first create an extended probability measure by
enriching the process through (carefully) choosing at random one of the particles to be the so called
\emph{spine}.
Now, on this enriched process,
changes of measure can easily be applied that \emph{only} affect the behaviour
along the path of this single distinguished `spine' particle;
in our case, we add a drift to the spine's motion, increase rate of fission along the path of the spine
and size-bias the spine's offspring distribution.
However, projecting this new enriched and changed measure back onto the original process filtration
(that is, without any knowledge of the distinguished spine)
brings the fundamental `additive' martingales into play as a Radon-Nikodym derivative.
The four probability measures, various martingales, extra filtrations and clear process constructions
afforded by our setup, together with some other useful properties and tricks, such as the
\emph{spine decomposition}, provides a very elegant, intuitive and powerful set of techniques for
analysing the process.

All three models that we consider in this article shall be built
on the same underlying space of sample trees with spines, and the
measures will all be constructed in analogous ways. Here we lay
out the details for the current model of branching Brownian
motion, and leave it to the reader to bridge the details when it
comes to the other models -- \cite{spine_foundations}
contains all the details in a more abstract setting that will
cover every model considered in this article.

The set of Ulam-Harris labels is to be equated with the set
$\Omega$ of finite sequences of strictly-positive integers:
\[
\Omega := \set{\emptyset} \cup \bigcup_{n\in\N} (\N)^n,
\]
where we take $\N=\set{1,2,\ldots}$. For two words $u,v\in\Omega$,
$uv$ denotes the concatenated word ($u\emptyset = \emptyset u =
u)$, and therefore $\Omega$ contains elements like `$213$' (or
`$\emptyset 213$'), which we read as `the individual being the 3rd
child of the 1st child of the 2nd child of the initial ancestor
$\emptyset$'. For two labels $v, u \in \Omega$ the notation $v <
u$ means that $v$ is an \emph{ancestor} of $u$, and $\abs{u}$
denotes the length of $u$. The set of all ancestors of $u$ is
equally given by
\[
\set{v:v<u} = \set{v: \exists w \in \Omega \text{ such that
}vw=u}.
\]

Collections of labels, ie. subsets of $\Omega$, will therefore be
groups of individuals. In particular, a subset $\tau\subset\Omega$
will be called a \emph{Galton-Watson tree} if:
\begin{enumerate}
\item $\emptyset\in\tau$,

\item if $u,v\in\Omega$, then $uv\in\tau$ implies $u\in\tau$,

\item for all $u\in\tau$, there exists $A_u\in{0,1,2,\ldots}$ such
that $uj\in\tau$ if and only if $1\leq j \leq 1+A_u$, (where
$j\in\N$).
\end{enumerate}
In general, each node $u$ of the underlying trees will have $1 + A_u$ branches coming from it,
whereas the special case of binary-branching, for example, would correspond
to $A_u = 1$ at every node.

The set of all Galton-Watson trees will be called $\mathbb{T}$.
Typically we use the name $\tau$ for a particular tree, and
whenever possible we will use the letters $u$ or $v$ or $w$ to
refer to the labels in $\tau$, which we may also refer to as
\emph{nodes of $\tau$} or \emph{individuals in $\tau$} or just as
\emph{particles}.

\bigskip\noindent
Each individual should have a \emph{location} in $\Real$ at each
moment of its \emph{lifetime}. Since a Galton-Watson tree
$\tau\in\mathbb{T}$ in itself can express only the \emph{family}
structure of the individuals in our branching process, in
order to give them these extra features we suppose that each
individual $u\in\tau$ has a mark $(X_u,\sigma_u)$ associated with
it which we read as:
\begin{itemize}
\item $\sigma_u\in\Real^+$ is the \emph{lifetime} of $u$, which
determines the \emph{fission time} of particle $u$ as $S_u :=
\sum_{v \leq u} \sigma_v$ (with $S_\emptyset:= \sigma_\emptyset$).
The times $S_u$ may also be referred to as the \emph{death} times;

\item $X_u:[S_u-\sigma_u,S_u)\to \Real$ gives the \emph{location}
of $u$ at time $t\in[S_u - \sigma_u,S_u)$.
\end{itemize}
To avoid ambiguity, it is always necessary to decide whether a
particle is in existence or not at its death time.
\begin{rem}\label{rem:death_and_birth_time_convention}
Our convention throughout will be that a particle $u$ dies `just
before' its death time $S_u$ (which explains why we have defined
$X_u:[S_u - \sigma_u,S_u)\to \cdot$ for example). Thus at the time
$S_u$ the particle $u$ has \emph{disappeared}, replaced by its $1+A_u$
children which are all alive and ready to go.
\end{rem}

We denote a single marked tree by $(\tau,X,\sigma)$ or $(\tau, M)$
for shorthand, and the set of all marked Galton-Watson trees by
$\mathcal{T}$:
\begin{itemize}
\item $\mathcal{T} := \bbset{(\tau,X,\sigma):
\tau\in\mathbb{T}\text{ and for each }u\in\tau,
\sigma_u\in\Real^+, X_u:[S_u - \sigma_u,S_u)\to \Real}$.

\item For each $(\tau, X, \sigma) \in \mathcal{T}$, the set of
particles that are alive at time $t$ is defined as $N_t :=
\set{u\in\tau: S_u - \sigma_u \leq t < S_u}$.
\end{itemize}

\bigskip\noindent
For any given marked tree $(\tau, M) \in \mathcal{T}$ we can
identify distinguished lines of descent from the initial ancestor:
$\emptyset, u_1, u_2, u_3,\ldots \in \tau$, in which $u_3$ is a
child of $u_2$, which itself is a child of $u_1$ which is a child
of the original ancestor $\emptyset$. We'll call such a subset of
$\tau$ a \emph{spine}, and will refer to it as $\xi$:
\begin{itemize}
\item a spine $\xi$ is a subset of nodes $\set{\emptyset, u_1,
u_2, u_3,\ldots}$ in the tree $\tau$ that make up a unique line of
descent. We use $\xi_t$ to refer to the unique node in $\xi$ that
that is alive at time $t$.
\end{itemize}
In a more formal definition, which can for example be found in the
paper by Rouault and Liu
\cite{rouault_liu:two_measures_on_tree_boundary}, a spine is
thought of as a point on $\partial\tau$ the boundary of the tree
-- in fact the boundary is \emph{defined} as the set of all
infinite lines of descent. This explains the notation $\xi \in
\partial\tau$ in the following definition: we augment the space
$\mathcal{T}$ of marked trees to become
\begin{itemize}
\item $\tilde{\mathcal{T}} := \bset{(\tau,M, \xi):
(\tau,M)\in\mathcal{T}\text{ and }\xi\in \partial\tau}$ is the set
of \emph{marked trees with distinguished spines}.
\end{itemize}

\bigskip\noindent
It is natural to speak of the \emph{position of the spine at time
$t$} which think of just as the position of the unique node that
is in the spine and alive at time $t$:
\begin{itemize}
\item we define the time-$t$ position of the spine as $\xi_t :=
X_u(t)$, where $u\in \xi \cap N_t$.
\end{itemize}
By using the notation $\xi_t$ to refer to both the node in the
tree and that node's spatial position we are introducing potential
ambiguity, but in practice the context will make clear which we
intend. However, in case of needing to emphasize, we shall give
the node a longer name:
\begin{itemize}
\item $\node{t}((\tau,M,\xi)) := u$ if $u \in \xi$ is the node in
the spine alive at time $t$,
\end{itemize}
which may also be written as $\node{t}(\xi)$.

As the spine $\xi_t$ diffuses, at the fission times $S_u$ for $u
\in \xi$ it gives birth to some offspring, one of which continues
the spine whilst the others go off to create subtrees like copies
of the BBM. These times on the spine are especially important for
the later spine decomposition of the martingale $Z_\lambda$, and
we therefore give them a name:
\begin{itemize}
\item the sequence of random times $\set{S_u: u \in \xi}$ are
known as the \emph{fission times on the spine};
\end{itemize}
Finally, it will later be important to know how many fission times
there have been in the spine, or what is the same, to know which
generation of the family tree the node $\xi_t$ is in (where the
original ancestor $\emptyset$ is considered to be the 0th
generation)
\begin{defn}\label{defn:fission_counting}
We define the counting function
\[
n_t = \abs{\node{t}(\xi)},
\]
which tells us which generation the spine node is in, or
equivalently how many fission times there have been on the spine.
For example, if $\xi_t = \seq{\emptyset, u_1, u_2}$ then both
$\emptyset$ and $u_1$ have died and so $n_t = 2$.
\end{defn}

The collection of all marked trees with a distinguished spine
is $\tilde{\mathcal{T}}$; on
this space we define four filtrations of key importance that encapsulate different
knowledge, but see Hardy and Harris \cite{spine_foundations} for more precise details:
\begin{itemize}
\item $\F_t$ knows everything that has happened to all the
branching particles up to the time $t$, \emph{but does not know
which one is the spine};

\item $\tilde{\F}_t$ knows everything that $\F_t$ knows and also
knows which line of descent is the spine (it is in fact the finest
filtration);

\item $\mathcal{G}_t$ knows only about the spine's motion in $J$
up to time $t$, but does not actually know which line of descent
in the family tree makes up the spine;

\item $\tilde{\mathcal{G}}_t$ knows about the spine's motion and
also knows which nodes it is composed of. Furthermore it knows
about the fission times of these nodes and how many children were
born at each time.
\end{itemize}

\bigskip\noindent
Having now defined the underlying space for our probabilities, we
remind ourselves of the probability measures:
\begin{defn}
For each $x \in \Real$, let $P^x$ be the measure on
$(\tilde{\mathcal{T}},\F_\infty)$ such that the filtered
probability space $(\tilde{\mathcal{T}},\F_\infty,(\F_t)_{t \geq
0}, P^\cdot)$ makes the $\Real$-valued point process $\mathbb{X}_t
= \set{X_u(t): u \in N_t}$ the canonical model for BBM.
\end{defn}

Our spine approach relies first on building a measure $\tilde{P}^x$
under which the spine is a single genealogical line of descent
chosen from the underlying tree. If we are given a
sample tree $(\tau,M)$ for the branching process,
it is easy to verify that,
if at each fission we make a uniform choice amongst the offspring
to decide which line of descent continues the spine $\xi$,
when $u \in \tau$ we have
\begin{equation}\label{uniform_spine_choice_definition}
\text{Prob}(u \in \xi) = \prod_{v < u} \frac{1}{1 + A_v}.
\end{equation}
This simple observation at
\eqref{uniform_spine_choice_definition} is the key to our method
for extending the measures, and for this we make use of the
following representation found in Lyons
\cite{lyons:simple_biggins_convergence}.
\begin{thm}
If $f$ is a $\tilde{\F}_t$-measurable function then we can write:
\begin{equation}\label{lyons_f_tilde_representation}
f = \sum_{u\in N_t} f_u \Ind_{(\xi_t = u)}
\end{equation}
where $f_u$ is $\F_t$-measurable.
\end{thm}
We use this representation to extend the measures $P^x$.
\begin{defn}\label{defn_of_P_tilde}
Given the measure $P^x$ on $(\tilde{\mathcal{T}}, \F_\infty)$ we
extend it to the probability measure $\tilde{P}^x$ on
$(\tilde{\mathcal{T}},\tilde{\F}_\infty)$ by defining
\begin{equation}\label{abcd}
\int_{\tilde{\mathcal{T}}} f \,\,\diffd \tilde{P}^x :=
\int_{\tilde{\mathcal{T}}} \sum_{u \in N_t} f_u \prod_{v < u} \frac{1}{1 +
A_v} \,\,\diffd P^x,
\end{equation}
for each $f\in m\tilde{\F}_t$ with representation like
\eqref{lyons_f_tilde_representation}.
\end{defn}
The previous approach to spines, exemplified in Lyons
\cite{lyons:simple_biggins_convergence}, used the idea of
\emph{fibres} to get a measure analogous to our $\tilde{P}$ that
could measure the spine. However, a perceived weakness in this approach was
that the corresponding measure had a time dependent total mass and
could not be normalized to become a probability measure
with an intuitive construction, unlike our $\tilde{P}$.
Our new idea of using the down-weighting
term of \eqref{uniform_spine_choice_definition} in the definition
of $\tilde{P}$ is crucial in ensuring that we get a
probability measure, and leads to the very useful situation in
which \emph{all} measure changes in our formulation are carried
out by \emph{martingales}.
\begin{thm}\label{p-tilde_is_extension}
This measure $\tilde{P}^x$ really is an extension of $P^x$ in that
$P = \tilde{P}\vert_{\F_\infty}$.
\end{thm}

\bigbreak

\noindent
The spine diffusion $\xi_t$ is $\tilde{\F}_t$-measurable, and it
is immediate that, under $\tilde{P}^x$, the spine diffusion $\xi_t$ is a Brownian
motion that starts at $x$.
In fact, it is easy to see that
\begin{thm}
Under $\tilde{P}^x$, the process $\mathbb{X}_t$ evolves as follows:
\begin{itemize}
\item starting from $x$, the spine $\xi_t$ diffuses
according to a (driftless) Brownian motion on $\Real$;

\item at rate $r$ the spine undergoes
fission producing $1+{A}$ particles, where ${A}$ is independent of the spine's motion
with distribution $\{ p_k:k\geq 0\}$;

\item with equal probability, one of the spine's offspring particles is
selected to continue the path of the spine,
repeating stochastically the behaviour of its parent;

\item the other particles initiate, from their birth position,
independent copies of a $P^\cdot$ branching Brownian motion with
branching rate $r$ and family-size distribution also given by $A$.
\end{itemize}
\end{thm}
It is useful to have this natural construction of $\tilde{P}$ in mind
when we start to change measure.

\subsection{New measures for BBM}
Having seen the construction of the underlying space and the
measure $\tilde{P}^x$, we can define a measure
$\tilde{\Q}_\lambda$ via a Radon-Nikodym derivative with respect
to $\tilde{P}$ after recalling three simpler changes of measures.

Changing measure with the $\tilde{P}$-martingale
\[
e^{\lambda (\xi_t-x) - \half \lambda^2 t}
\]
would make the spine process $\xi_t$ a Brownian motion with drift $\lambda$ under the new measure.

We also note that
\[
e^{-mrt}(1+m)^{n_t}
\]
is a $\tilde{P}$-martingale that will increase the rate of the Poisson process $n_t$ of fission
times on the spine from $r$ to $(1+m)r$ after a change of measure.

Lastly,
\[
\prod_{v<\xi_t}\left(\frac{1+A_v}{1+m}\right)
\]
is also a $\tilde{P}$-martingale that will produce size-biasing (only) along the spine,
giving an offspring distribution $\{{(1+k)}p_k /{(1+m)}: k\geq 0\}$.

Combining these components into a single change of measure leads to the following key result:
\begin{thm}
\label{thm:tildeQ_definition}
Define the measure $\tilde{\Q}_\lambda$ on
$(\tilde{\mathcal{T}},\tilde{\F}_\infty)$ by:
\begin{equation}
\label{highbury}
\frac{\diffd \tilde{\Q}_\lambda^x}{\diffd\tilde{P}^x}\bigg\vert_{\tilde{\F}_t}
=
e^{\lambda (\xi_t-x) - E_\lambda t}
\,
\prod_{v<\xi_t}\left({1+A_v}\right).
\end{equation}
Under $\tilde{\Q}_\lambda^x$, the process $\mathbb{X}_t$ evolves as follows:
\begin{itemize}
\item starting from $x$, the spine $\xi_t$ diffuses
according to a Brownian motion with drift $\lambda$ on $\Real$;

\item at accelerated rate $(1+m)r$ the spine undergoes
fission producing $1+\tilde{A}$ particles, where $\tilde{A}$
is independent of the spine's motion
with size-biased distribution $\{ (1+k)p_k/(1+m):k\geq 0\}$;

\item with equal probability, one of the spine's offspring particles is
selected to continue the path of the spine,
repeating stochastically the behaviour of its parent;

\item the other particles initiate, from their birth position,
independent copies of a $P^\cdot$ branching Brownian motion with
branching rate $r$ and family-size distribution given by $A$, that is,
$\{p_k:k\geq 0\}$.
\end{itemize}
\end{thm}

The measure $\Q_\lambda$ that we introduced in Theorem \ref{thm:BBM_c-o-m_with_Z}
via its pathwise construction can equivalently be obtained from
$\tilde{\Q}_\lambda$ by restricting it to the filtration
$(\F_t)_{t \geq 0}$ on which the original measure $P$ is defined, that is,
simply ignoring information identifying the spine:
\begin{thm}\label{thm:BBM_c-o-m_with_Z_version_tilde}
If we define $\Q_\lambda^x :=
\tilde{\Q}_\lambda^x\vert_{{\F}_\infty}$, then $\Q_\lambda^x$ is a
measure on $\F_\infty$ with
\begin{equation}
\frac{\diffd \Q_\lambda^x}{\diffd P^x}\bigg\vert_{\F_t} =
\frac{Z_\lambda(t)}{Z_\lambda(0)}.
\label{Qcofm}
\end{equation}
Under ${\Q}_\lambda^x$,
the branching-diffusion point process
$\mathbb{X}_t$ can be pathwise constructed exactly as described in
Theorem \ref{thm:BBM_c-o-m_with_Z}.
\end{thm}
There are at least two ways to prove this result: Kyprianou
\cite{kyprianou:travelling_waves} bases his proof on a
decomposition of a related (non-probability) measure ${P}^*$
as a product of measures for the spine's motion, the fission-counting process $n_t$, and
measures on the sub-trees born from the spine.
Since we have clear constrictions of the probability measures $\tilde{P}$ and $\tilde{\Q}_\lambda$
and richer collection of filtrations, we have an alternative:

\bigskip\noindent\textbf{Proof of Theorem
\ref{thm:BBM_c-o-m_with_Z_version_tilde}: }
Since $\Q_\lambda^x := \tilde{\Q}_\lambda^x\vert_{{\F}_\infty}$,
whilst $P^x =\tilde{P}^x\vert_{{\F}_\infty}$ from Theorem \ref{p-tilde_is_extension}
and it is clear that the change of
measure at \eqref{highbury} projects onto the sub-algebra $\F_t$ as a
conditional expectation, we have
\[
\frac{\diffd \Q_\lambda^x}{\diffd P^x}\bigg\vert_{\F_t}
=
\frac{\diffd \tilde{\Q}_\lambda^x}{\diffd
\tilde{P}^x}\bigg\vert_{\F_t} = e^{-\lambda x} \,
\tilde{P}^x
\ccond{e^{\lambda \xi_t - E_\lambda t}
\,\prod_{v<\xi_t}\left({1+A_v}\right)}{\F_t}.
\]
Using representation \eqref{lyons_f_tilde_representation}
and recalling from \eqref{uniform_spine_choice_definition}
that $\tilde{P}\cond{\xi_t = u}{\F_t} = \prod_{v < u} \left({1+A_v}\right)^{-1}$
yields
\begin{align*}
\tilde{P}^x\ccond{ e^{\lambda \xi_t - E_\lambda t}\,\prod_{v<\xi_t}\left({1+A_v}\right)}{\F_t}
&=
\tilde{P}^x\ccond{\sum_{u \in N_t}
e^{\lambda X_u(t) - E_\lambda t} \times \prod_{v<u}\left({1+A_v}\right) \times \Ind_{(\xi_t = u)}}{\F_t}
\\
&= \sum_{u \in N_t} e^{\lambda X_u(t) - E_\lambda t} \times \prod_{v<u}\left({1+A_v}\right)
 \times \tilde{P}\cond{\xi_t = u}{\F_t}
\\
&=  \sum_{u \in N_t} e^{\lambda X_u(t) - E_\lambda t}
= Z_\lambda(t)
\end{align*}
and trivially noting $Z_\lambda(0)=e^{\lambda x}$ under $P^x$ completes the proof.
\qed

\subsection{Proof of Theorem \ref{BBM_L_p_theorem}}\label{sec:proof_of_BBM}

Just before we proceed to the proof we recall the naturally
occurring eigenvalue $E_\lambda := \half \lambda^2 + mr$,
noting that
under the symmetry assumption that $\lambda\leq 0$
and for $p\in(1,2]$:
\[
p E_\lambda - E_{p\lambda} > 0 \quad\Leftrightarrow\quad
c_\lambda > c_{p\lambda}
\quad\Leftrightarrow\quad p \lambda^2 < 2mr
\]
and that this always holds for some $p>1$ whenever $\lambda\in (\tilde\lambda,0]$,
that is, when $\lambda$ lies between the minimum of $c_\lambda$ found at $\tilde\lambda$ and the origin.

\subsubsection{Proof of part 1:} We are going to prove that for
every $p\in(1,2]$ the martingale $Z_\lambda$ is
$\eL^p(P)$-convergent if $p E_\lambda - E_{p\lambda} > 0$.
Furthermore, since $P^x(Z_\lambda(t)^p) = e^{p\lambda x}
P^0(Z_\lambda(t)^p)$ we do not lose generality supposing that $x=0$;
from now on this is implicit
if we drop the superscript by simply writing $P$.

From the change of measure in Theorem \ref{thm:BBM_c-o-m_with_Z}
or Theorem \ref{thm:BBM_c-o-m_with_Z_version_tilde} it is clear
that
\[
P(Z_\lambda(t)^p) = P(Z_\lambda(t)^{p-1} Z_\lambda(t)) =
\Q_\lambda (Z_\lambda(t)^q),
\]
where $q:=p-1$. Our aim is to prove that $\Q_\lambda
(Z_\lambda(t)^q)$ is bounded in $t$, since then
$Z_\lambda^p(t)$ must be bounded in $\eL^p(P)$ and Doob's theorem
will then imply that $Z_\lambda$ is convergent in $\eL^p(P)$.

As we mentioned, the algebra $\tilde{\mathcal{G}}_\infty$ gives us
the very important \emph{spine-decomposition} of the martingale
$Z_\lambda$:
\begin{equation}
\tilde{\Q}_\lambda\cond{Z_\lambda(t)}{\tilde{\mathcal{G}}_\infty}
= \sum_{k=1}^{n_t} A_k e^{\lambda \xi_{S_k}- E_\lambda S_k} +
e^{\lambda \xi_t-E_\lambda t},\label{spine_decomp}
\end{equation}
where $A_k$ is the number of new particles produced from the fission at time
$S_k$ along the path of the spine, and the sum is taken to equal 0 if $n_t = 0$.
Full details of this can be found in \cite{spine_foundations}, but the
intuition is quite clear:
since the particles that do not make up the spine
grow to become independent copies of $\mathbb{X}_t$
distributed \emph{as if under} $P$, the fact that $Z_\lambda$ is a
$P$-martingale on these subtrees implies that their contributions
to the above decomposition are just equal to their
\emph{immediate} contribution on being born at time $S_k$ at
location $\xi_{S_k}$. \Remark We emphasize that here we must use
the measure $\tilde{\Q}_\lambda$, since $\Q_\lambda$ cannot
measure the algebra $\tilde{\mathcal{G}}_\infty \nsubseteq
\F_\infty$.

We can now use the conditional form of Jensen's inequality
followed by the spine decomposition of \eqref{spine_decomp}
coupled with the simple inequality,
\begin{prop}\label{neveu_pythagoras_inequality}
If $q\in(0,1]$ and $u,v>0$ then $(u+v)^q\leq u^q + v^q$,
\end{prop}
to obtain,
\begin{align}
\tilde{\Q}_\lambda\cond{Z_\lambda(t)^q}{\tilde{\mathcal{G}}_\infty}
& \leq
\tilde{\Q}_\lambda\cond{Z_\lambda(t)}{\tilde{\mathcal{G}}_\infty}^q
\label{cond_jensen_with_q}
\\
& \leq
\sum_{k=1}^{n_t} A_k^q  e^{q\lambda \xi_{S_k}- qE_\lambda S_k} +
e^{q\lambda \xi_t-qE_\lambda t}.
\end{align}
With the tower property of conditional expectations and noting that $\Q_\lambda$
and $\tilde{\Q}_\lambda$ agree on $\mathcal{F}_t$,
\begin{align}
\Q_\lambda(Z_\lambda(t)^q) = \tilde{\Q}_\lambda(Z_\lambda(t)^q)
& =
\tilde{\Q}_\lambda\left(\tilde{\Q}_\lambda\cond{Z_\lambda(t)^q}{\tilde{\mathcal{G}}_\infty}\right)
\\
&\leq \tilde{\Q}_\lambda \bseq{\sum_{k=1}^{n_t} A_k^q  e^{q\lambda
\xi_{S_k}- qE_\lambda S_k}} + \tilde{\Q}_\lambda \bseq{e^{q\lambda
\xi_t-qE_\lambda t}},\label{here_i_am}
\end{align}
and the proof of $\eL^p(P)$-boundedness will be complete once we show
this is bounded in $t$.

As written, \eqref{here_i_am} is made up of two terms, and since
they play a central role from here on we name them explicitly: on
the far right we have the \textbf{spine term} $\tilde{\Q}_\lambda
\seq{e^{q\lambda \xi_t-qE_\lambda t}}$, the other being the
\textbf{sum term} $\tilde{\Q}_\lambda \seq{\sum_{k=1}^{n_t}
A_k^q  e^{q\lambda \xi_{S_k}- qE_\lambda S_k}}$.

\bigskip\noindent\textbf{The spine term:}
Changing from $\tilde{P}$ to $\tilde{\Q}_\lambda$ gives the spine
a drift of $\lambda$, and therefore the change-of-measure for just
the spine's motion (\ie on the algebra $\mathcal{G}_t$) is carried
out by the martingale $e^{\lambda \xi_t - \half\lambda^2 t}$, so
\begin{align}
\tilde{\Q}_\lambda \bseq{e^{q\lambda \xi_t-qE_\lambda t}} &=
\tilde{P}
\bseq{e^{q\lambda \xi_t-qE_\lambda t} \times e^{\lambda \xi_t - \half\lambda^2 t}}\notag\\
&= e^{\{\half(p\lambda)^2 -\half\lambda^2 \}t - qE_\lambda t }\, \tilde{P}
\bseq{e^{p\lambda \xi_t - \half(p\lambda)^2 t}}\notag\\
&= e^{ -(p E_\lambda - E_{p\lambda})t}\,\tilde\Q_{p\lambda}(1)
=e^{ -(p E_\lambda - E_{p\lambda})t}\label{right_here_right_now}
\end{align}
since the second-line term $e^{p\lambda \xi_t - \half(p\lambda)^2 t}$
is also a $\tilde{P}$-martingale
and $\half(p\lambda)^2-\half\lambda^2=E_{p\lambda}-E_\lambda$.

\bigskip\noindent\textbf{The sum term:} Recall, under the measure $\tilde{\Q}_\lambda$
we know that the fission times $\{S_k:k\geq0\}$ on the spine occur as a
Poisson process of rate $(1+m)r$ with the $k^{\mathrm{th}}$ fission yielding an additional $A_k$ offspring,
each $A_k$ being an \emph{independent} copy of $\tilde{A}$ which has the size-biased
distribution $\{(1+k)p_k/(1+m):k\geq0\}$.

First, conditioning on the motion of the spine (without knowledge of the fission times or
family sizes) and appealing to intuitive results from
Poisson process theory (see \cite{kallenberg} for example) yields
\begin{equation}
\tilde{\Q}_\lambda
\ccond{\sum_{k=1}^{n_t} A_k^q e^{q\lambda \xi_{S_k}- qE_\lambda S_k}}
{\,\mathcal{G}_t}
=
 \int_0^t (1+m)r \,
 \tilde{\Q}_\lambda\seq{\tilde{A}^q} \,
 e^{q\lambda \xi_s- qE_\lambda s} \, \diffd s
\label{poppy}
\end{equation}
Taking expectations of both sides of \eqref{poppy}
and using Fubini's theorem then gives
\begin{align}
\tilde{\Q}_\lambda\bseq{ \sum_{k=1}^{n_t} A_k^q e^{q\lambda
\xi_{S_k}- qE_\lambda S_k}}
&= (1+m)r\,\tilde{\Q}_\lambda\seq{\tilde{A}^q}
\int_0^t  \, \tilde{\Q}_\lambda \bseq{
e^{q\lambda
\xi_s- qE_\lambda s}} \, \diffd s \notag\\
&= (1+m)r \, \tilde{\Q}_\lambda(\tilde{A}^q) \int_0^t e^{-(p E_\lambda
- E_{p\lambda})s} \, \diffd s, & \text{using
}\eqref{right_here_right_now}.\notag
\end{align}
\bigskip\noindent Thus we have found an explicit upper-bound
(if $pE_\lambda \neq E_{p\lambda}$):
\begin{equation}
P^x(Z_\lambda(t)^p) \leq e^{p\lambda x}
\left(\frac{(1+m)r}{pE_\lambda - E_{p\lambda}}
\bsseq{1 - e^{-(p E_\lambda - E_{p\lambda})t}}\tilde{\Q}_\lambda\seq{\tilde{A}^q} + e^{-(p E_\lambda
- E_{p\lambda})t}\right).
\end{equation}
Finally, we also observe that
\begin{lem}\label{BBM_family_size_l_p}
If $p\in(1,2]$ and $q:=p-1$, $\tilde{\Q}_\lambda(\tilde{A}^q)<\infty$ if and only if $P(A^p)<\infty$
\end{lem}
since
\[
\tilde{\Q}_\lambda(\tilde{A}^q) = \sum_{i=1}^\infty i^q
\frac{i+1}{m+1} p_i = \frac{1}{m+1} \bseq{\sum_{i=1}^\infty i^p
p_i + \sum_{i=1}^\infty i^q p_i} = \frac{P(A^p) + P(A^q)}{m+1} \leq
\frac{2 P(A^p)}{m+1}.
\]

Hence, if we have $pE_\lambda - E_{p\lambda} > 0$ in addition to $P(A^p)<\infty$,
this implies that $P^x(Z_\lambda(t)^p)$ will remain bounded as $t \to \infty$, which
together with Doob's theorem will complete the proof of the first
part of Theorem \ref{BBM_L_p_theorem}. \qed

\subsubsection{Proof of Part 2:}

We seek to show that $Z_\lambda$ is unbounded in $\mathcal{L}^p(P^x)$ if
either $pE_\lambda-E_{p\lambda}<0$ or $P(A^p)=\infty$.

Note that \emph{if} $Z_\lambda$ is $\mathcal{L}^p(P^x)$ bounded then
\[
P^x(Z_\lambda(\infty)^p)=\lim_{t\rightarrow\infty}P^x(Z_\lambda(t)^p)<\infty
\]
hence, $\tilde{\Q}_\lambda^x(Z_\lambda(\infty)^q)<\infty$ and $Z_\lambda(\infty)^q$
is a uniformly integrable $\tilde{\Q}_\lambda^x$-submartingale.
In particular, for any stopping time $T$,
$\tilde{\Q}_\lambda^x\cond{Z_\lambda(\infty)^q}{\mathcal{F}_T}\geq Z_\lambda(T)^q$
hence
$\tilde{\Q}_\lambda^x(Z_\lambda(\infty)^q)\geq \tilde{\Q}_\lambda^x(Z_\lambda(T)^q)$.

First, by considering only the contribution of the spine
$Z_\lambda(t)\geq e^{\lambda\xi_t-E_\lambda t}$ for all $t\geq 0$ and
recalling \eqref{right_here_right_now}, we see that
\[
\tilde{\Q}_\lambda^x(Z_\lambda(t)^q)\geq \tilde{\Q}_\lambda^x(e^{q\lambda\xi_t-qE_\lambda t})
= e^{q\lambda x-(pE_\lambda - E_{p\lambda})t}
\]
and $Z_\lambda$ is therefore unbounded in $\mathcal{L}^p(P^x)$ if $pE_\lambda-E_{p\lambda}<0$.

Now, let $T$ be any fission time along the path of the spine, then
\[
Z_\lambda(T)\geq (1+\tilde{A})e^{\lambda \xi_T - E_\lambda T}
\]
where $\tilde{A}$ is the number of additional offspring produced at the time of fission. Then,
\[
\tilde{\Q}_\lambda^x(Z_\lambda(T)^q)\geq
\tilde{\Q}_\lambda^x\left((1+\tilde{A})^q\right)\,e^{q\lambda x}\,\tilde{\Q}_{p\lambda}^x(e^{-(pE_\lambda-E_{p\lambda})T})
\]
and so $Z_\lambda$ is unbounded in $\mathcal{L}^p(P^x)$ if
$\tilde{\Q}_\lambda^x\left((1+\tilde{A})^q\right)=\infty\iff P(\tilde{A}^p)=\infty$.
\qed

\section{A typed branching diffusion}\label{chapter_MTBBM}

We move on to consider a \emph{typed} branching diffusion where the type of each particle
evolves as a Markov chain and influences the fission rate, offspring distribution
and spatial diffusion coefficient of the particle.
We will follow a similar notation and setup as for the BBM case, but leave some details to
the reader.

\bigskip\noindent
A \textbf{`typical' particle motion} 
will be like the process $(X_t,Y_t)_{t\geq0}$ in $J:=\mathbb{R}\times I$ where: \\
\ (i) the \emph{type} location, $Y_t$, is an irreducible,
time-reversible Markov chain on the finite type-space $I:=\{1,\ldots,n\}$ with Q-matrix $\theta Q$
($\theta$ is a strictly positive constant) and invariant measure
$\pi=(\pi_1,\ldots,\pi_n)$; \\
\ (ii) the \emph{spatial} location, $X_t$, moves as
a driftless Brownian motion on $\Real$ with diffusion coefficient $a(y)>0$ whenever
$Y_t$ is in state $y$, that is,
\begin{equation}
\diffd X_t = a(Y_t)^\half \,\diffd B_t,\qquad \textrm{ where }B_t\textrm{ a Brownian motion}.
\end{equation}
The formal generator of this process $(X_t,Y_t)$ is therefore:
\begin{equation}\label{eqn:SPM_generator}
\mathcal{H}F(x,y)=\half a(y)\frac{\partial^2 F}{\partial x^2} +
\theta \sum_{j\in I} Q(y,j) F(x,j),\qquad(F:J\to\Real).
\end{equation}
We often use matrix calculations, and it is convenient to gather
the diffusion coefficients together in a diagonal matrix
$A:=\textrm{diag}[a(1),\ldots,a(n)]$.

\bigskip\noindent\textbf{The typed branching diffusion.}
Consider a typed branching Brownian motion whereby individual particles move independently
according to the \emph{`typical' particle motion}, as described above,
and at any moment a particle currently of type $y$ will undergo fission
at rate $r(y)$ to be replaced by a random number of offspring, $1+A(y)$,
where $A(y)\in \set{0,1,2,\ldots}$ is an independent RV with distribution
\[
P\seq{A(y) = i} = p_i(y), \qquad i \in \set{0,1,\ldots},
\]
and mean $m(y):=P\seq{A(y)}<\infty$ for all $y\in I$.
At birth, offspring inherit the parent's spatial and type positions and then move off independently,
 repeating stochastically the parent's behaviour, and so on.
As before, we let $N(t)$ be the set of particles alive at time $t$ using the Ulam-Harris labelling,
where a particle $u\in N(t)$ has spatial position and type given by $(X_u(t),Y_u(t))$.
We gather together the birth rates in a diagonal matrix
$R:=\textrm{diag}[r(1),\ldots,r(n)]$ and the mean number of offspring in
$M:=\textrm{diag}[m(1),\ldots,m(n)]$.

We suppose that the probabilities for this process are given by $\set{P^{x,y}: (x,y) \in J}$
defined on the natural filtration $(\F_t)_{t \geq 0}$, where $P^{x,y}$ is the law of the process
starting with one initial particle of type $y$ at spatial position $x$.

\bigskip\noindent\textbf{The typed branching diffusion with spine.}
As seen before, we can extend the measures $\set{P^{x,y}: (x,y) \in J}$ to
$\set{\tilde{P}^{x,y}: (x,y) \in J}$ by identifying a single distinguished infinite line of
descent starting from the initial single particle, known as the \emph{spine}.
The natural filtration $(\tilde{\F}_t)_{t \geq 0}$ then contains all information
about the process, including the identity of the spine.
 \begin{thm}\label{thm:MT_original_spine_construction}
Under $\tilde{P}^{x,y}$, the process $\mathbb{X}_t:=\{(X_u(t),Y_u(t)) : t\geq0\}$ evolves as follows:
\begin{itemize}
\item starting from $(x,y)$, the \emph{spine} $(\xi_t,\eta_t)$ evolves as
a process with type $\eta_t$ moving as a Markov chain on $I$ with $Q$-matrix $\theta Q$
and spatial position $\xi_t$ diffusing as a (driftless) Brownian on $\Real$
with variance coefficient $a(\eta_t)$;

\item if the spine is currently of type $y$, at rate $r(y)$ the spine undergoes
fission producing $1+{A(y)}$ particles, where ${A(y)}$ is an independent random variable
with distribution $\{ p_k(y):k\geq 0\}$;

\item with equal probability, one of the spine's offspring particles is
selected to continue the path of the spine,
repeating stochastically the behaviour of its parent;

\item the other particles initiate, from their birth position,
independent copies of a $P^{\cdot,\cdot}$ typed branching Brownian motion
as described above.
\end{itemize}
\end{thm}

\bigskip
It should be noted that the condition of time-reversibility on the
Markov chain is not absolutely necessary, and is really just a
simplifying assumption that gives us an easier $\eL^2$ theory for
the matrices and eigenvectors;
our aim is really to show how the spine techniques work
-- lessening the geometric complexity of the model serves a good purpose.

Note, the special case of the 2-type BBM model was considered
in Champneys \emph{et al} \cite{art:aap} by different means.
Also, in our model, at the time of fission a type-$y$
individual can produce only type-$y$ offspring. This is not the
same as the case in which a type-$y$ individual may produce a
random collection of particles of \emph{different} types --
as considered in T.E.\ Harris's
classic text \cite{teharris}, for example.
Other forms of typed branching processes have also been dealt with
by spine techniques, for example,
see Lyons \emph{et al} \cite{lyons_kurtz_pemantle_peres:conceptual_Kesten_Stigum}
or
Athreya \cite{athreya:bernoulli_change_of_measure_MC} for discrete-time
models in which a particle's type does not change during its life
but a type-$w$ individual can give offspring of any type according
to some distribution.
See also the remarkable work of Georgii and Baake \cite{georgii} that uses spine techniques
to study ancestral type behaviour in a  continuous time branching Markov chain
where particles can give birth to across all types.
In principle, our spine methods will be robust enough to extend to all these other type behaviours
(with added spatial diffusion).

\subsection{The martingale}

Via the many-to-one lemma (see \cite{spine_foundations}) or generators
it is easy to see that for any $\lambda \in
\Real$, any function (vector) $v_\lambda:I \to \Real$ and any
number $E_\lambda \in \Real$, the expression
\begin{equation}
Z_\lambda(t) := \sum_{u\in N(t)} v_\lambda(Y_k(t)) \, e^{\lambda
X_k(t)-E_\lambda t},
\end{equation}
will be a martingale if and only if $v_\lambda$ and $E_\lambda$
satisfy:
\begin{equation}
\bseq{\half\lambda^2A+\theta Q+MR} v_\lambda = E_\lambda
v_\lambda,\label{eigenvector_condition1}
\end{equation}
that is, $v_\lambda$ must be an eigenvector of the
matrix $\half\lambda^2A+\theta Q+MR$, with eigenvalue $E_\lambda$.
\begin{defn}
For two vectors $u,v$ on $I$, we define
\[
\ip{u}{v}_\pi := \sum_{i=1}^n u_i v_i \pi_i,
\]
which gives us a Hilbert space which we refer to as $\eL^2(\pi)$.
We suppose that the eigenvector $v_\lambda$ is normalized so that
$\norm{v_\lambda}_\pi := \ip{v_\lambda}{v_\lambda}_\pi = 1$.
\end{defn}
The fact that the Markov chain is time-reversible implies that the
matrix $\half\lambda^2A+\theta Q+MR$ is self-adjoint with respect
to this inner product. This in itself is enough to guarantee the
existence of eigenvectors in $\eL^2(\pi)$, but the fact that we
are dealing with a finite-state Markov chain means that we also
have the \emph{Perron-Frobenius} theory to hand, which allows us
to suppose that $v_\lambda$ is a \emph{strictly positive}
eigenvector whose eigenvalue $E_\lambda$ is real and the farthest
to the right of all the other eigenvalues -- see Seneta
\cite{seneta} for details. This implies a useful representation
for the eigenvalue:
\begin{thm}
\begin{equation}\label{sup_rep_of_eigenvalue}
E_\lambda=\sup_{\norm{v}_\pi=1} \ip{\left((\lambda^2/2)A+\theta Q
+MR\right)v}{v}_\pi,
\end{equation}
since it is the rightmost eigenvalue.
\end{thm}
A proof can be found in Kreyzig \cite{kreyszig}. From this it is
not difficult to show that $E_\lambda$ is a strictly-convex
function of $\lambda$.
Interestingly, it will be seen
in our proofs that it is the geometry of the eigenvalue
$E_\lambda$ that determines the interval that gives rise to
martingales $Z_\lambda(t)$ that are $\eL^p$-convergent.
\begin{cor}\label{e_lambda_facts}
As a function of $\lambda$, $E_\lambda$ is strictly-convex and
infinitely differentiable with
\begin{equation}
E_\lambda' = \lambda \ip{A v_\lambda}{v_\lambda}_\pi.
\end{equation}
If we define the speed function
\begin{equation}\label{c_lambda_definition}
c_\lambda := -E_\lambda/\lambda,
\end{equation}
then on $(-\infty,0)$ the function $c_\lambda$ has just one
minimum at a single point $\tilde{\lambda}(\theta)$,
either side of which $c_\lambda$ is strictly increasing to $+\infty$ as
either $\lambda\downarrow-\infty$ or $\lambda\uparrow0$.
In particular,
for each $\lambda \in\insidewith$ there is some $p> 1$ such that
$c_\lambda>c_{p\lambda}$;
on the other hand,
if $\lambda\outsidewithout$ there is no such $p > 1$.
\end{cor}
We refer to the function $c_\lambda$ as the speed function since
it relates to the asymptotic speed of the travelling waves
associated with the martingale $Z_\lambda(t)$; see Harris
\cite{art:twFKPP} or Champneys \emph{et al} \cite{art:aap} for
details of the relationship between branching-diffusion
martingales and travelling waves.

Since $Z_\lambda(t)$ is a strictly-positive martingale it is
immediate that $Z_\lambda(\infty):= \lim_{t \to \infty}
Z_\lambda(t)$ exists and is finite almost-surely under $P^{x,y}$.
As before, \emph{by symmetry we shall assume that $\lambda\leq 0$ and,
without loss of generality, we also suppose that $P(A(y)=0)=1$ whenever $r(y)=0$
to simplify statements}.
We shall prove necessary and sufficient conditions for
$\eL^1$-convergence of the $ Z_\lambda$ martingales:

\begin{thm}\label{MTBBM_cgce_theorem}
For each $x\in \Real$, the limit $Z_\lambda(\infty) := \lim_{t \to \infty} Z_\lambda(t)$ exists
$P^{x,y}$-a.s. where:
\begin{itemize}
\item  if $\lambda \outsidewith$ then $Z_\lambda(\infty) = 0$ $P^{x,y}$-almost surely;
\label{MTkyp_part_1}
\item  if $\lambda \in\insidewith$ and $P(A(y) \log^+A(y))=\infty$ for some $y\in I$,
then $Z_\lambda(\infty) = 0$ $P^{x,y}$-a.s.;
\label{MTkyp_part_2}
\item if $\lambda \in\insidewith$ and $P(A(y) \log^+A(y))<\infty$ for all $y\in I$,
then $Z_\lambda(t)\rightarrow Z_\lambda(\infty)$ almost surely and in $\eL^1(P^{x,y})$.
\label{MTkyp_part 3}
\end{itemize}
\end{thm}

Once again, in many cases where the martingale has a non-trivial limit,
the convergence will be much stronger than merely in $\eL^1(P^{x,y})$, as
indicated by the following new $\eL^p$-convergence result that we will prove by extending our
earlier new spine approach:
\begin{thm}\label{thm:finite_type_convergence}
For each $x \in \Real$, and for each $p \in (1, 2]$:
\begin{itemize}
\item $Z_\lambda(t)\rightarrow Z_\lambda(\infty)$ a.s. and in $\eL^p(P^{x,y})$
if $p E_\lambda - E_{p\lambda} > 0$ and $P(A(y)^p)<\infty$ for all $y\in I$.
\label{MTLp_part_1}
\item $Z_\lambda$ is unbounded in $\eL^p(P^{x,y})$, that is, $P^{x,y}(Z_\lambda(\infty)^p)=\infty$ if
either $p E_\lambda - E_{p\lambda} < 0$ or $P(A(y)^p)=\infty$ for some $y\in I$.
\label{MTLp_part_2}
\end{itemize}
Note, when $\lambda\leq 0$, the inequality $p E_\lambda - E_{p\lambda} > 0$
is equivalent to $c_\lambda>c_{p\lambda}$ and holds for some $p\in(1,2]$ if and only if
$\lambda\in\insidewith$.
\end{thm}

\subsection{New measures for the typed BBM}
With the construction of the underlying space and the
measure $\tilde{P}^{x,y}$,  we can define a measure
$\tilde{\Q}_\lambda$ via a Radon-Nikodym derivative with respect
to $\tilde{P}$ by combining three simpler changes of measures that only affect
behaviour along the spine.

We first alter the motion of the spine:
\begin{lem}
\[
v_\lambda(\eta_t) \, e^{\int_0^t MR(\eta_s) \, \diffd s}
e^{\lambda \xi_t -E_\lambda t}
\]
is a $\tilde{P}$-martingale and
using this martingale as a Radon-Nikodym derivative would give a new measure
under which $\xi_t$ would have instantaneous drift $a(\eta_t)\lambda$
and $\eta_t$ would have modified $Q$-matrix $Q_\lambda$
with invariant measure $\pi_\lambda = v_\lambda^2 \pi$.
\end{lem}

It is easy to see this first martingale property with classical `one-particle' calculations,
for example, using the Feynman-Kac formula and noting the relation \eqref{eigenvector_condition1}.
We will discuss the behaviour of the spine under such a change of measure in detail in the next subsection.

Secondly, we alter the fission rate along the spine:
\begin{lem}
\[
e^{-\int_0^t MR(\eta_s) \, \diffd s} \prod_{v < \xi_t}
\seq{1+m(\eta_{S_v})}
\]
is a $\tilde{P}$-martingale that will increase the rate at which fission
times occur on the spine from $R(\eta_t)$ to
$(1+m(\eta_t))R(\eta_t)$.
\end{lem}

Finally, we size-bias the offspring distribution along the spine:
\begin{lem}
\[
\prod_{v < \xi_t} \frac{1+A_v}{1+m(\eta_{S_v})}
\]
is a $\tilde{P}$-martingale that will cause the family distribution on the
spine to be size-biased to the distribution
\[
\text{Prob}(\tilde{A}(w) = i) = \frac{(i+1)p_i(w)}{m(w)+1},\qquad i \in
\set{0,1,\ldots}.
\]
\end{lem}

Combining these components in one change of measure leads to the following key definition:
\begin{defn}\label{MTBBMCofM}
For each $\lambda \in \Real$ we define a measure
$\tilde{\Q}_\lambda^{x,y}$ on $(\tilde{\mathcal{T}},
\tilde{\F}_\infty)$ via
\begin{equation}
\frac{\diffd \tilde{\Q}_\lambda^{x,y}}{\diffd
\tilde{P}^{x,y}}\bigg\vert_{\tilde{\F}_t}
:=
\frac{v_\lambda(\eta_t)}{v_\lambda(y)}
\, {e^{\lambda (\xi_t-x)-E_\lambda t}}
  \, \prod_{v < \xi_t} (1 + A_v).
\end{equation}
\end{defn}
As in the previous BBM case, a proof using the conditional expectation of this
measure-change martingale confirms:
\begin{thm} \label{MTBBM_QCofM}
If we define $\Q_\lambda^{x,y} :=
\tilde{\Q}_\lambda^{x,y}\vert_{{\F}_\infty}$, then
$\Q_\lambda^{x,y}$ is a measure on $\F_\infty$ and
\[
\frac{\diffd \Q_\lambda^{x,y}}{\diffd P^{x,y}}\bigg\vert_{\F_t} =
\frac{Z_\lambda(t)}{Z_\lambda(0)}
\]
\end{thm}

Notice that with this result, starting with the three simple known martingales,
we have actually shown that $Z_\lambda$ must, in fact, be a martingale.
This route offers a simple way of getting very general `additive' martingales for
the branching process.

\subsection{The spine process $(\xi_t, \eta_t)$ under $\tilde{\Q}_\lambda$}

In the BBM model it was clear to see that the spine
$\xi_t$ received a drift under the measure $\tilde{\Q_\lambda}$; something
similar happens here:

\begin{lem}\label{the_facts}
Under $\tilde{\Q}_\lambda$ the spine process $(\xi_t, \eta_t)$ has
generator:
\begin{equation}\label{eqn:SPM_generator_lambda}
\mathcal{H}_\lambda F \, (x,y):=\half a(y)\frac{\partial^2
F}{\partial x^2} + a(y) \lambda \frac{\partial F}{\partial x} +
\sum_{j\in I} \theta Q_\lambda(y,j) F(x,j),
\end{equation}
where $Q_\lambda$ is an honest Q-matrix:
\[
\theta Q_\lambda(i,j) = \left\{\begin{array}{cc} \theta Q(i,j)
\frac{v_\lambda(j)}{v_\lambda(i)} & \textrm{if }i\neq j\\
\theta Q(i,i)+\frac{\lambda^2}{2}a(i) - E_\lambda + r(i) &
\textrm{if }i=j
\end{array}\right.
\]
with invariant measure $\pi_\lambda = v_\lambda^2 \pi$.
\end{lem}
Thus under $\tilde{\Q}_\lambda$ the Q-matrix (generator) of
$\eta_t$ is changed to $\theta Q_\lambda$, and the process $\xi_t \in \Real$ is given an
instantaneous drift of
$a(\eta_t) \lambda$.
The form of this above generator can be obtained from the theory
of Doob's $h$-transforms, due to the fact that on the algebra
$\mathcal{G}_t$ the change of measure is given by:
\begin{equation}\label{finite_type_spine_c-o-m}
\frac{\diffd \Q_\lambda^{x,y}}{\diffd
P^{x,y}}\bigg\vert_{\mathcal{G}_t}
=
\frac{1}{v_\lambda(y)e^{\lambda x}}
\, v_\lambda(\eta_t) \, e^{\int_0^t MR(\eta_s) \,
\diffd s} e^{\lambda \xi_t -E_\lambda t}.
\end{equation}

The long-term behaviour under $\tilde{\Q}_\lambda$ of the spine
diffusion $\xi_t$ can now be retrieved from the generator
\eqref{eqn:SPM_generator_lambda} and the properties of $E_\lambda$
stated in Lemma \ref{e_lambda_facts}:
\begin{cor}\label{drift_thm}
Almost surely under $\tilde{\Q}_\lambda^{x,y}$, the long-term drift of the spine is given explicitly as
\[
\lim_{t\to\infty} t^{-1} \xi_t = E_\lambda'
\]
and hence
\begin{equation}
\xi_t+c_\lambda t\rightarrow
\begin{cases}
\infty & \textrm{if } \lambda\in\insidewith \\
-\infty & \textrm{if } \lambda\outsidewithout
\end{cases}
\end{equation}
whereas, if $\lambda=\tilde{\lambda}$ the process $\xi_t +c_\lambda t$ will be recurrent on $\Real$
under $\tilde\Q_\lambda$.
\end{cor}
\Proof From the generator stated at
\eqref{eqn:SPM_generator_lambda} we can write:
\[
\xi_t = B\left(\int_0^t a(\eta_s) \diffd s\right) + \lambda
\int_0^t a(\eta_s) \diffd s,
\]
where $B(t)$ is a $\tilde{\Q}_\lambda$-Brownian motion. Then by
the ergodic theorem and the fact that $\pi_\lambda = v_\lambda^2
\pi$:
\[
t^{-1} \xi_t \to \lambda \sum_{y\in I} a(y) \pi_\lambda(y) =
\lambda \sum_{y\in I} a(y) v_\lambda^2(y) \pi(y) = \lambda \ip{A
v_\lambda}{v_\lambda}_\pi = E_\lambda'.
\]
Direct calculation from \eqref{c_lambda_definition} gives
$E_\lambda' = -c_\lambda - \lambda c_\lambda'$, and therefore
$t^{-1} (\xi_t + c_\lambda t) \to -\lambda c_\lambda'$, whence
whether we are to the left or right of the local minimum of $c_\lambda$
found at $\tilde\lambda$ determines the behaviour of $\xi_t + c_\lambda t$,
as is required. Lastly, when $\lambda=\tilde{\lambda}$,
with the laws of the iterated logarithm in mind,
it is not difficult to see that  both $B(\int_0^t a(\eta_s)\,ds)$ and
$\int_0^t (\lambda a(\eta_s) + c_{{\lambda}})\,ds$ will fluctuate about the origin,
hence $\xi_t + c_{\lambda} t$
will be recurrent under $\tilde\Q_\lambda$.
\qed

\subsection{Construction of the process under $\tilde{\Q}_\lambda$}

Drawing together the elements from this section, we now present the important intuitive
pathwise construction of the new measure $\tilde{\Q}_\lambda^{x,y}$:

\begin{thm}
\label{thm:tildeQ_construction}
Under $\tilde{\Q}_\lambda^{x,y}$, the process $\mathbb{X}_t$ evolves as follows:
\begin{itemize}
\item starting from $(x,y)$, the spine $(\xi_t,\eta_t)$ evolves as a Markov process with
generator $\mathcal{H}_\lambda$,
that is, $\eta_t$ evolves as Markov chain on $I$ with $Q$-matrix $\theta Q_\lambda$
and $\xi_t$ moves as a Brownian motion on $\Real$ with variance coefficient $a(\eta_t)$
and drift $a(\eta_t)\lambda$.
\item whenever the type of the spine $\eta$ is in state $y\in I$,
the spine undergoes fission at an accelerated rate $(1+m(y))r$,
producing $1+\tilde{A}(y)$ particles where $\tilde{A}(y)$
is independent of the spine's motion with size-biased distribution
$\{ (1+k)p_k(y)/(1+m(y)):k\geq 0\}$;
\item with equal probability, one of the spine's offspring particles is
selected to continue the path of the spine,
repeating stochastically the behaviour of its parent;
\item the other particles initiate, from their birth position,
independent copies of $P^{\cdot,\cdot}$ typed branching Brownian motions.
\end{itemize}
\end{thm}

\subsection{Proof of Theorem \ref{MTBBM_cgce_theorem}}
The following proof is an extension of that given for BBM by Kyprianou
\cite{kyprianou:travelling_waves}.
The second part of the following theorem is
the key element in using the measure change \eqref{MTBBMCofM}
to determine properties of the martingale $Z_\lambda$:
\begin{thm}\label{measure_result}
Suppose that $P$ and $\Q$ are two probability measures on a space
$\seq{\Omega, \F_\infty}$ with filtration $(\F_t)_{t\geq 0}$, such
that for some positive martingale $Z_t$,
\[
\frac{\diffd \Q}{\diffd P}\bigg\vert_{\F_t} = Z_t.
\]
The limit $Z_\infty := \limsup_{t \to \infty} Z_t$ therefore
exists and is finite almost surely under $P$. Furthermore, for any
$F \in \F_\infty$
\begin{equation}\label{durret_decomposition}
\Q(F) = \int_F Z_\infty \, \diffd P + \Q\seq{F \cap \{Z_\infty =
\infty\}},
\end{equation}
and consequently
\begin{align}
&(a) \quad P(Z_\infty =0)=1  \iff \Q(Z_\infty = \infty)=1 \label{zeromeasure}\\
&(b) \quad P(Z_\infty)= 1 \iff \Q(Z_\infty < \infty)=1 \label{UImeasure}
\end{align}
\end{thm}
A proof of the decomposition \eqref{durret_decomposition} can be
found in Durrett \cite{durrett:prob_with_examples}, at page 241.

Suppose that $\lambda\leq \tilde{\lambda}<0$.
Ignoring all contributions except for the spine, it is immediate that
\[
Z_\lambda(t)
= \sum_{u \in N_t} v_\lambda(Y_u(t))\, e^{\lambda X_u(t)- E_\lambda t}
\geq v_\lambda(\eta_t)\, e^{\lambda(\xi_t + c_\lambda t)}
\]
where, from Corollary \ref{drift_thm}, under the measure $\tilde{\Q}_\lambda$ the spine
satisfies $\liminf\{\xi_t+c_\lambda t\}= -\infty$ a.s. and $v_\lambda>0$, hence
$\limsup_{t \to \infty}Z_\lambda(t) =\infty$ almost surely under $\tilde{\Q}_\lambda$,
yielding $P(Z_\lambda(\infty)=0)=1$.

Note that, for $y\in I$,
$P(A(y)\log^+{A})<\infty \iff \sum_{k\geq 1} P(\log^+{\tilde{A}(y)}>ck)<\infty$ for any $c>0$,
where recall that $\tilde{A}(y)$ has the size-biased distribution $\{(i+1)p_k(y)/(1+m(y)):k\geq0\}$.
Then for an IID sequence $\{\tilde{A}_n(y)\}$ of copies of $\tilde{A}(y)$, Borel-Cantelli reveals that,
$P$ almost surely,
\begin{equation}
\limsup_{t\rightarrow\infty} n^{-1}\log^+{\tilde{A}_n(y)}
=
\begin{cases}
0 & \text{if } P(A(y) \log^+{A}(y))<\infty, \\
\infty & \text{if } P(A(y)\log^+{A}(y))=\infty.\\
\end{cases}\label{expon_A_behaviour}
\end{equation}

Now suppose that $\lambda\in(\tilde{\lambda},0]$ and $P(A(y)\log^+{A}(y))=\infty$
for some $y\in I$ (with $r(y)>0$).
Let $S_k$ be the time of the $k^{\textrm{th}}$ fission along the spine
producing $\tilde{A_k}(\eta_{S_k})$ additional particles, then
\[
Z_\lambda(S_k)\geq \tilde{A}_{k}(\eta_{S_k}) v_\lambda(\eta_{S_k})\,e^{\lambda(\xi_{S_k}+c_\lambda S_k)}
\]
where $(\xi_t + c_\lambda t)/t\rightarrow -\lambda c_\lambda^\prime>0$,
$\eta_t$ is ergodic so the event $\{\eta_{S_k}=y\}$ will occur for infinitely many $k$
since $r(y)>0$,
and $n_t/t\rightarrow <Rv_\lambda,v_\lambda>_\pi$ so $S_k/k \rightarrow  <Rv_\lambda,v_\lambda>_\pi^{-1}$,
hence the super-exponential growth for $\tilde{A}_k(y)$ from \eqref{expon_A_behaviour} gives
$\limsup_{t \to \infty}Z_\lambda(t) =\infty$ $\tilde{\Q}_\lambda$-almost surely
which then implies that $P(Z_\lambda(\infty)=0)=1$.

Finally, suppose that $\lambda\in(\tilde{\lambda},0]$ and $P(A(y)\log^+{A}(y))<\infty$
for all $y\in I$.
Recall the very important \emph{spine-decomposition} of the martingale
$Z_\lambda$ from \eqref{spine_decomp}:
\begin{equation}
\tilde{\Q}_\lambda\cond{Z_\lambda(t)}{\tilde{\mathcal{G}}_\infty}
= \sum_{k=1}^{n_t} \tilde{A}_{k}(\eta_{S_k}) v_\lambda(\eta_{S_k})\,e^{\lambda(\xi_{S_k}+c_\lambda S_k)} +
v_\lambda(\eta_{t})e^{\lambda (\xi_t+c_\lambda t)}.
\end{equation}
In this case, 
the facts that $(\xi_t + c_\lambda t)/t\rightarrow -\lambda c_\lambda^\prime>0$
and $S_k/k \rightarrow  <Rv_\lambda,v_\lambda>_\pi^{-1}$ together with
the moment conditions and \eqref{expon_A_behaviour} implying that the $\tilde{A}_k(y)$'s
all have sub-exponential growth
means that
\[
\limsup_{t\rightarrow\infty} \tilde{\Q}_\lambda\cond{Z_\lambda(t)}{\tilde{\mathcal{G}}_\infty}
<\infty \qquad \tilde\Q_\lambda \textrm{-a.s.}
\]
Fatou's lemma then gives $\liminf_{t\rightarrow\infty} Z_\lambda(t)<\infty$, $\tilde\Q_\lambda$-a.s., hence also
$\Q_\lambda$-a.s. In addition, since $Z_\lambda(t)^{-1}$ is a positive $\Q_\lambda$-martingale
(recall Theorem \ref{MTBBM_QCofM}) with an almost sure limit, this means that $\lim_{t\rightarrow\infty} Z_\lambda(t)<\infty$, $\Q_\lambda$-a.s.
 and then \eqref{UImeasure} yields that $P(Z_\lambda(\infty))=1$ and so $Z_\lambda(t)$ converges almost surely and in
$\mathcal{L}^1(P)$.
\qed

\subsubsection{Discussion of rate of convergence to zero and left-most particle speed.}
Alternatively, when $\lambda<\tilde{\lambda}$ we
can readily obtain the rate of convergence to zero with the following simple argument,
adapted from Git \emph{et al} \cite{ghh}.
By Proposition \ref{neveu_pythagoras_inequality},
\[
Z_\lambda(t)^q
\leq \sum_{u\in N(t)} v_\lambda(Y_u(t))^q\, e^{q\lambda(X_u(t)+c_{q\lambda})} \,
e^{ q\lambda(c_\lambda-c_{q\lambda})t}
\leq K\, Z_{q\lambda}(t) e^{q\lambda(c_\lambda-c_{q\lambda})t}
\]
where $K:=\max_{y\in I} v_\lambda^q(y)/v_{q\lambda}(y)<\infty$ since $I$ is finite and
$v_\lambda>0$.
Recall that $c_\lambda$ has a minimum over $\lambda\in(-\infty,0]$ at $\tilde\lambda$
with $c_{\tilde\lambda}=-E_{\tilde\lambda}=-\tilde{\lambda}<Av_{\tilde\lambda},v_{\tilde\lambda}>_\pi$.
Then, since $Z_{q\lambda}(t)$ is a convergent martingale,
we can choose $q$ such that $q\lambda=\tilde{\lambda}$ giving $Z_\lambda(t)$ decaying exponentially
to zero at least at rate $\lambda(c_\lambda-c_{\tilde\lambda})$.

Further, once we know that $P$ and $\Q_\lambda$ are equivalent for every $\lambda\in(\tilde\lambda,0]$,
since the spine moves such that $\xi_t/t\rightarrow -c_\lambda-\lambda c_\lambda^{\prime}$ under $\Q_\lambda$,
the left-most particle $L(t):=\inf_{u\in N(t)}X_u(t)$ must satisfy $\liminf_{t} L(t)/t \leq -c_{\tilde\lambda}$, $P$-a.s.
On the other hand, the convergence of the $Z_\lambda$ $P$-martingales quickly
gives the same upper bound on the fastest speed of any particle,
leading to $L(t)/t\rightarrow -c_{\tilde\lambda}$, $P$-a.s.
This result also reveals that the rate of exponential decay found above is actually best possible.

\subsection{Proof of Theorem \ref{thm:finite_type_convergence}}

\subsubsection{Proof of Part 1:} Suppose $p\in(1,2]$, then with
$q:=p-1$ a slight modification of the BBM proof arrives at
\begin{align*}
P^{x,y}(Z_\lambda(t)^p)
&=
e^{\lambda x}v_\lambda(y)\tilde\Q_\lambda^{x,y}(Z_\lambda(t)^q)\\
&\leq e^{\lambda x}v_\lambda(y)\left(
\tilde\Q_\lambda^{x,y}
\bseq{\sum_{k=1}^{n_t} A_k^q
v_\lambda(\eta_{S_k})^q e^{q\lambda \xi_{S_k}- qE_\lambda S_k}} +
\tilde\Q_\lambda^{x,y} \bseq{v_\lambda(\eta_t)^q e^{q\lambda
\xi_t-qE_\lambda t}}\right)
\end{align*}
and the proof of $\eL^p$-boundedness will be complete once we show
that this RHS expectation is bounded in $t$.

\bigskip\noindent\textbf{The spine term.} Since the type space $I$ is finite, we trivially note that
$\ip{v_\lambda^p}{v_{p\lambda}}_\pi<\infty$.
 It is always useful to first focus on
the spine term, since we can change the measure with
\eqref{finite_type_spine_c-o-m} to get
\begin{align}
\tilde{\Q}_\lambda^{x,y} \bseq{v_\lambda(\eta_t)^q e^{q\lambda
\xi_t-qE_\lambda t}} &= \tilde{P}^{x,y} \bbseq{v_\lambda(\eta_t)^q
e^{q\lambda \xi_t-qE_\lambda t} \times \frac{ v_\lambda(\eta_t) \,
e^{\int_0^t M R(\eta_s) \, \diffd s} e^{\lambda \xi_t -E_\lambda
t}}{v_\lambda(y)
e^{\lambda x}}}\notag\\
&= e^{q\lambda x} \frac{v_{p\lambda}(y)}{v_\lambda(y)} g_t(y)
e^{-(pE_\lambda - E_{p\lambda}) t}
\label{spine_term_multitype_BBM}
\end{align}
where, for all $y\in I$
\[
g_t(y):=\tilde{\Q}_{p\lambda}^{0,y}
\bseq{\frac{v_\lambda^p}{v_{p\lambda}}(\eta_t)}
\rightarrow \ip{v_{\lambda}^p}{v_{p\lambda}}_\pi
\]
as $t\rightarrow\infty$ and
$\ip{g_t v_{p\lambda}}{v_{p\lambda}}_\pi=\ip{v_{\lambda}^p}{v_{p\lambda}}_\pi$ for all $t\geq0$,
since $\eta_t$ is a finite-state irreducible Markov
chain under $\tilde\Q_{\mu}$ with invariant distribution
$\pi_\mu(y) = v_\mu(y)^2 \pi(y)$.
It follows
that the long term the growth or decay of the spine term is
determined by the sign of $pE_\lambda - E_{p\lambda}$.

\bigskip\noindent\textbf{The sum term.}
We now assume that $p E_\lambda-E_{p\lambda}>0$.
We know that under $\tilde{\Q}_\lambda$ and conditional on knowing $\eta$,
the fission times $\{S_k:k\geq 0\}$ on the spine occur
as a Poisson process of rate $(1+m(\eta_s))r(\eta_s)$
with the $k^{\textrm{th}}$ fission yielding an additional
$A_k$ offspring, each $A_k$ being an independent copy of $\tilde{A}(y)$
which has the size-biased distribution $\{ (1+k)p_k(y)/(1+m(y)) : k\geq 0\}$
where $y=\eta_{S_k}$ is the type at the time of fission.
We also recall from Lemma \ref{BBM_family_size_l_p} that
\[
M_q(y):=\tilde{\Q}_\lambda(\tilde{A}^q(y))<\infty \iff P(A^p(y))<\infty.
\]

Therefore, if we condition on $\mathcal{G}_t$ which
knows about $(\xi_s,\eta_s)$ at all times $0 \leq s \leq t$ we can
transform the sum into an integral, use Fubini's theorem
and the change of measure used in \eqref{spine_term_multitype_BBM}:
\begin{align*}
\tilde\Q_\lambda^{x,y}
\bseq{\sum_{k=1}^{n_t} A_k^q v_\lambda(\eta_{S_k})^q
e^{q\lambda \xi_{S_k}- qE_\lambda S_k}}
&= \tilde\Q_\lambda^{x,y}
\bseq{\tilde\Q_\lambda^{x,y}\ccond{\sum_{k=1}^{n_t} A_k^q v_\lambda(\eta_{S_k})^q
e^{q\lambda \xi_{S_k}- qE_\lambda S_k}}{\mathcal{G}_t}}\notag \\
&= \tilde\Q_\lambda^{x,y} \bseq{\int_0^t (1+m(\eta_s)) r(\eta_s) \,
M_q(\eta_s) v_\lambda(\eta_s)^q
e^{q\lambda \xi_s- qE_\lambda s} \, \diffd s}\notag \\
&= \int_0^t \tilde\Q_\lambda^{x,y} \bseq{(1+m(\eta_s)) r(\eta_s)
M_q(\eta_s) {v_\lambda(\eta_s)^q e^{q\lambda \xi_s- qE_\lambda s}}} \, \diffd
s \\
&=
e^{q\lambda x} \frac{v_{p\lambda}(y)}{v_\lambda(y)} \int_0^t
h_s(\eta_s)
e^{-(pE_\lambda - E_{p\lambda}) s} \, \diffd s
\\
&=e^{q\lambda x} \frac{v_{p\lambda}(y)}{v_\lambda(y)}\times\,
\frac{k_t(y)}{pE_\lambda-E_{p\lambda}}
\end{align*}
where
\[
h_s(y):=
\tilde\Q_{p\lambda}^{0,y} \bseq{\tilde{r}(\eta_s) M_q(\eta_s) \,
\frac{v_\lambda^p}{v_{p\lambda}}(\eta_s)},
\qquad
\tilde{r}(y):=(1+m(y))r(y),
\]
\[
\textrm{and }\qquad k_t(y):=\Expct(h_U(y);U\leq t)
\]
with $U$ an independent exponential of rate $(pE_\lambda-E_{p\lambda})>0$.
Note that, for all $y\in I$, $h_s(y)\rightarrow \ip{\tilde{r} M_q v_\lambda^p}{v_{p\lambda}}$
and $k_t(y)\uparrow k_\infty(y)$ as $t\rightarrow\infty$, where
$\ip{k_t v_{p\lambda}}{v_{p\lambda}}=\ip{\tilde{r} M_q v_\lambda^p}{v_{p\lambda}}\,\Prob(U\leq t)
\uparrow \ip{k_\infty v_{p\lambda}}{v_{p\lambda}}=\ip{\tilde{r} M_q v_\lambda^p}{v_{p\lambda}}$.
Then, since $M_q(w)<\infty \iff P(A(w)^p)<\infty$, and $I$ is finite, we are guaranteed that
$k_\infty(y)<\infty$ for all $y\in I$ as long as $P(A(w)^p)<\infty$ for all $w\in I$.

\bigskip\noindent Having dealt with both the spine term and the sum term, we have
obtained the upper-bound
\[
P^{x,y}\seq{Z_\lambda(t)^p} \leq
\frac{e^{p\lambda x}v_{p\lambda}(y)}{\left(p E_\lambda-E_{p\lambda}\right)}
\left(
k_t(y) + g_t(y) \left(p E_\lambda-E_{p\lambda}\right) e^{-(p E_\lambda-E_{p\lambda})t}
\right)
\]
and since $Z_\lambda(t)^p$ is a $P$-submartingale, we find that
\[
P^{x,y}\seq{Z_\lambda(t)^p} \leq
\frac{e^{p\lambda x}v_{p\lambda}(y)}{\left(p E_\lambda-E_{p\lambda}\right)}
\,k_\infty(y)
\qquad (\forall t\geq 0)
\]
and $Z_\lambda(t)$ will be bounded in $\eL^p(P^{x,y})$
if we have both $pE_\lambda - E_{p\lambda}>0$ and $P(A^p(w))<\infty$ for all $w\in I$.
\qed

\subsubsection{Proof of Part 2:}
The earlier proof for BBM goes through with minor modification.
Exactly as in the BBM case,
looking only at the contribution of the spine means that
$Z_\lambda$ is unbounded in $\mathcal{L}^p(P^{x,y})$
if $pE_\lambda-E_{p\lambda}<0$.
In addition, letting $T$ be any fission time along the path of the spine,
\[
Z_\lambda(T)\geq (1+\tilde{A}(\eta_T))v_\lambda(\eta_T)e^{\lambda \xi_T - E_\lambda T}
\]
where $\tilde{A}(\eta_T)$ is the number of additional offspring produced at the time of fission.
Then, with $m_q(y):=\tilde{\Q}((1+\tilde{A}(y))^q)<\infty \iff P(A^p(y))<\infty$,
\begin{align*}
\tilde{\Q}_\lambda^{x,y}(Z_\lambda(T)^q)
&\geq
e^{q\lambda x}\,\tilde{\Q}_\lambda^{x,y}(m_q(\eta_T) v_\lambda(\eta_T)^q e^{q\lambda \xi_T-qE_\lambda T}) \\
& = e^{q\lambda x}\,\tilde{\Q}_{p\lambda}^{x,y}\left(m_q(\eta_T) \frac{v_\lambda(\eta_T)^p}{v_{p\lambda}(\eta_T)}
\, e^{-(pE_\lambda-E_{p\lambda}) T}\right)
\end{align*}
and so $Z_\lambda$ will also be unbounded in $\mathcal{L}^p(P^{x,y})$ if
$m_q(y)=\infty \iff P(A^p(y))=\infty$ for any $y\in I$ (taking a fission time when also in state $y$).
\qed

\subsubsection{Remarks on signed martingales and Kesten-Stigum type theorems}
In the multi-typed BBM, for each $\lambda$
there will be other (signed) additive martingales
corresponding to the different eigenvectors and eigenvalues obtained from solving
\eqref{eigenvector_condition1};
the $Z_\lambda$ martingale simply corresponds to the Perron-Frobenius, or ground-state,
eigenvalue $E_\lambda$ and (strictly positive) eigenvector $v_\lambda$.
Since $\abs{u+v}^q\leq (\abs{u}+\abs{v})^q\leq  \abs{u}^q+\abs{v}^q$ for all $u,v\in\Real$,
the above proof will also adapt to give convergence results for signed martingales.
In fact, when there is a complete orthonormal set of eigenvectors, a Kesten-Stigum
like theorem would then swiftly follow
(for example, see Harris \cite{harris:gibbs_boltzmann} in the context of the continuous-type
model of the next section).

\section{A continuous-typed branching-diffusion}\label{chapter:HW_model} %

The previous finite-type model was originally inspired by the
model that we now turn to, originally laid out in Harris and
Williams \cite{art:largedevs1}. In this model the type moves on
the real line as an Orstein-Uhlenbeck process associated with the
generator
\[
Q_\theta := \frac{\theta}{2} \bseq{ \frac{\partial^2}{\partial
y^2} - y \frac{\partial}{\partial y}},\qquad\text{with }\theta
>0\text{ considered as the \emph{temperature}},
\]
which has the standard normal density as its invariant
distribution:
\[
\pi(y) := (2 \pi)^{-\half} e^{-\half y^2}.
\]
The spatial movement of a particle of type $y$ is a driftless
Brownian motion with instantaneous variance
\[
A(y) := ay^2, \qquad\text{for some fixed }a > 0,
\]
and fission of a particle of type $y$ occurs at a rate
\[
R(y) := ry^2 + \rho,\qquad\text{where }r,\rho >0 \text{ are
fixed,}
\]
to produce two particles at the same type-space location as the
parent (we consider only binary splitting). The model has very
different behaviour for low temperature values (\ie low $\theta$),
but most studies have considered the high temperature regime where
$\theta > 8r$. Also, the parameter $\lambda$ must be restricted to
an interval $(\lambda_{\text{min}},0)$ in order for some of the
model's parameters to remain in $\Real$, where
\[
\lambda_{\text{min}} := - \sqrt{\frac{\theta - 8r}{4a}}.
\]
Generally, \emph{unboundedness} in a model's rates is a serious
obstacle to classical proofs since they often depend on the
expectation semigroup of the branching process, and unbounded
rates tend to lead to unbounded \emph{eigenfunctions}. Here this
is the case, but the existence of a spectral theory for their
particular expectation operator allowed Harris and Williams to get
a sufficiently good bound in particular for a \emph{non-linear} term (see
Theorem 5.1 of \cite{art:largedevs1}), and therefore to prove
$\eL^p$-convergence of the martingale.
Other convergence results for various martingales
and weighted sums over particles for this model
also appear in Harris \cite{harris:gibbs_boltzmann},
again using more classical methods and requiring `non-linear' calculations.
The spine approach we again adopt here is both simple and more generic in nature;
requiring no such special `non-linear'  calculations,
it elegantly produces very good estimates that only involve
easy one-particle calculations.

We use the same notation as previously $\mathbb{X}_t =
\set{\seq{X_u(t), Y_u(t)}: u \in N_t}$ to denote the point process
of space-type locations in $\Real \times \Real$, and suppose that
the measures $\set{\tilde{P}^{x,y}: (x,y) \in \Real^2}$ on the
natural filtration with a spine $(\tilde{\F}_t)_{t \geq 0}$ are
such that the initial ancestor starts at $(x,y)$ and
$\seq{\mathbb{X}_t, (\xi_t, \eta_t)}$ becomes the above-described
branching diffusion with a spine.

\subsection{The measure change}

Although there are some significant differences, this model is
similar in flavour to our finite-type model.
There is a strictly-positive
martingale $Z_\lambda$ defined as
\[
Z_\lambda(t) := \sum_{u \in N_t} v_\lambda(Y_u(t)) e^{\lambda
X_u(t) - E_\lambda t}
\]
where $v_\lambda$ and $E_\lambda$ are the eigenvector and
eigenvalue associated with the self-adjoint (in $\eL^2(\pi)$)
operator:
\[
Q_\theta + \half\lambda^2 A(y) + R(y).
\]
The  eigenfunction $v_\lambda$ is normalizable against the
$\eL^2(\pi)$ norm, and can be found explicitly as
\[
v_\lambda(y) = e^{\psi_\lambda^- y^2}
\]
where
\[
\psi_\lambda^- := \frac{1}{4} - \frac{\mu_\lambda}{2\theta}, \quad
\mu_\lambda := \half \sqrt{\theta^2 - \theta(8r + 4 a \lambda^2)},
\]
are both positive for all $\lambda \in (\lambda_{\text{min}},0)$;
another important parameter is $\psi_\lambda^+ := \frac{1}{4} +
\frac{\mu_\lambda}{2\theta}$.
The eigenvalue $E_\lambda$ is then given by
\[
E_\lambda = \rho + \theta \psi_\lambda^-.
\]
We again define the speed function $c_\lambda :=
-E_\lambda/\lambda$, and $\tilde{\lambda}(\theta) < 0$ is the
unique point (on the negative axis) at which $c_\lambda$ hits
its minimum -- further details are given in Harris and Williams
\cite{art:largedevs1}.
We are going to use spines to prove the following result,
in which the critical case of $\lambda=\tilde\lambda$ and
the necessary conditions for $\eL^p(P)$-convergence are new results:
\begin{thm}\label{thm:OU_type_convergence}
Suppose that $\lambda \in (\lambda_{\emph{min}},0)$.
\begin{enumerate}
\item Let $p\in(1,2]$. The martingale $Z_\lambda$ is
$\eL^p(P)$-bounded if both $p E_\lambda - E_{p\lambda} > 0$ and
$p\psi_\lambda^- < \psi_{p\lambda}^+$.
In particular, for all $\lambda\in\insidewith$, $Z_\lambda$ is a uniformly-integrable martingale.

\item $Z_\lambda$ is unbounded in $\eL^p(P)$ if either $p E_\lambda - E_{p\lambda} < 0$
or $p\psi_\lambda^- > \psi_{p\lambda}^+$.

\item Almost surely under $P$, $Z_\lambda(\infty) = 0$ if
$\lambda\outsidewith$.
\end{enumerate}
\end{thm}
Once again, for each $\lambda \leq 0$ we define a measure
$\tilde{\Q}_\lambda^{x,y}$ on $(\tilde{\mathcal{T}},
\tilde{\F}_\infty)$ via
\begin{equation}
\frac{\diffd \tilde{\Q}_\lambda^{x,y}}{\diffd
\tilde{P}^{x,y}}\bigg\vert_{\tilde{\F}_t} := \frac{1}{v_\lambda(y)
e^{\lambda x}} 2^{n_t} v_\lambda(\eta_t) e^{\lambda \xi_t
-E_\lambda^- t},
\end{equation}
so that with $\Q_\lambda := \tilde{\Q}_\lambda\vert_{\F_\infty}$
we have
\[
\frac{\diffd \Q_\lambda^{x,y}}{\diffd P^{x,y}}\bigg\vert_{\F_t} =
\frac{Z_\lambda(t)}{Z_\lambda(0)} =
\frac{Z_\lambda(t)}{v_\lambda(y) e^{\lambda x}}.
\]
The facts are that under $\tilde{\Q}_\lambda$:
\begin{itemize}
\item the spine diffusion $\xi_t$ has instantaneous drift $a
\eta_t^2 \lambda$;

\item the type process $\eta_t$ has generator $\frac{\theta}{2}
\seq{\frac{\partial^2}{\partial y^2} - \frac{2\mu_\lambda}{
\theta} y\frac{\partial}{\partial y}}$ and an invariant probability measure
$\pi_{\lambda}:=\ip{v_\lambda}{v_\lambda}_\pi^{-1} v_\lambda^2\pi$,
corresponding to a normal distribution,  $N(0,\frac{\theta}{2 \mu_\lambda})$;

\item fission times on the spine occur at the accelerated rate of
$2 R(\eta_t)$;

\item all particles not in the spine behave as if under the
original measure $P$.
\end{itemize}

\subsection{Proof of Theorem \ref{thm:OU_type_convergence}}
\bigskip\noindent\textbf{Proof of Part 1:}
Suppose $p\in(1,2]$. Then using the spine decomposition with
Jensen's inequality and Proposition
\ref{neveu_pythagoras_inequality} we find,
\begin{gather*}\label{multitype(general)_spine_decomp}
P^{x,y}(Z_\lambda^-(t)^p)
\leq e^{\lambda x}v_\lambda(y)\left(
\tilde{\Q}_\lambda^{x,y}
\bseq{\sum_{u<\xi_t} v_\lambda(\eta_{S_u})^q e^{q\lambda
\xi_{S_u}- qE_\lambda S_u}} + \tilde\Q_\lambda^{x,y}
\bseq{v_\lambda(\eta_t)^q e^{q\lambda \xi_t-qE_\lambda t}}\right).
\end{gather*}
Assume that $pE_\lambda - E_{p\lambda} > 0$ and $p\psi_\lambda^- < \psi_{p\lambda}^+$.
As seen in Harris and Williams \cite{art:largedevs1},
we can do many calculations explicitly in this model, largely due to the fact that
under $\tilde{\Q}_{p\lambda}^{0,y}$
\[
\eta_s\sim \text{N}\left(e^{- \mu_{p\lambda}s} y, \frac{\theta(1 - e^{-2
\mu_{p\lambda} s})}{2 \mu_{p\lambda}}\right) \rightarrow \text{N}\left(0,\frac{\theta}{2 \mu_{p\lambda}}\right)
\]
 and the eigenfunctions $v_\lambda$ have such simple exponential form.
For example,
\begin{equation}\label{OU_birthrate_term}
\tilde{\Q}_{p\lambda}^{0,y} \bseq{
\frac{v_\lambda^p}{v_{p\lambda}}(\eta_s)}
=
\tilde{\Q}_{p\lambda}^{0,y}
\bseq{e^{(p\psi_\lambda^- - \psi_{p\lambda}^-) \eta_s^2}}
\end{equation}
can easily be seen to be finite and bounded for all $s\geq 0$ if and only if
$p\psi_\lambda^- - \psi_{p\lambda}^- -
\frac{\mu_{p\lambda}}{\theta} = p\psi_\lambda^- -
\psi_{p\lambda}^+< 0$, and just as readily calculated explicitly.

In fact, more `natural' conditions for $\eL^p$ convergence of the martingales would be that
\[
\ip{ R M_q v_{\lambda}^p}{v_{p\lambda}}_\pi<\infty, \qquad \ip{ v_{\lambda}^p}{v_{p\lambda}}_\pi<\infty,
\qquad \text{and}\qquad pE_\lambda-E_{p\lambda}<0,
\]
where $M_q(y):=\tilde{\Q}(\tilde{A}^q(y))$ with $\tilde{A}$ the size-biased
offspring distribution (here, binary splitting means $\tilde{A}(y)\equiv 1$),
 and we present arguments below that are more generic in nature, at least in terms of
adapting to other `suitably' ergodic type motions and random family sizes.
Note, the last condition above is related to the natural convexity of $E_\lambda$ and,
in our specific model, both integrability conditions are guaranteed
by $p\psi_\lambda^- -\psi_{p\lambda}^+< 0$.

\bigskip\noindent\textbf{The spine term.}
On the algebra $\mathcal{G}_t$ the change of measure takes the
form
\[
\frac{\diffd \tilde{\Q}_\lambda^{x,y}}{\diffd
\tilde{P}^{x,y}}\bigg\vert_{\mathcal{G}_t}
=
\frac{v_\lambda(\eta_s)}{v_\lambda(y)}\, \exp{\left(\int_0^t R(\eta_s) \, \diffd s
+ \lambda (\xi_t-x) - E_\lambda^- t\right)},
\]
which we can use on the spine term to arrive
at
\begin{equation}\label{spine_term_OUBBM}
f_t(x,y) := e^{\lambda x}v_\lambda(y)
\tilde{\Q}_\lambda^{x,y} \bseq{v_\lambda(\eta_t)^q e^{q\lambda\xi_t-qE_\lambda t}}
=
e^{p\lambda x} v_{p\lambda}(y) \,
g_t(y) e^{-(pE_\lambda - E_{p\lambda}) t}
\end{equation}
with
$g_t(y):=\tilde{\Q}_{p\lambda}^{0,y} \bseq{{v_\lambda^p}(\eta_t)/{v_{p\lambda}}(\eta_t)}$.
Under the assumption that $p\psi_\lambda^- < \psi_{p\lambda}^+$, it easy to check that
$\ip{ v_{\lambda}^p}{v_{p\lambda}}_\pi<\infty$,
that is ${v_\lambda^p}/{v_{p\lambda}}\in L^1(\pi_{p\lambda})$ from which it follows that
$g_t\in L^1(\pi_{p\lambda})$ for all $t\geq0$.
Since $\eta$ has equilibrium $\pi_{p\lambda}$ under $\tilde \Q_{p\lambda}$,
we find $\ip{g_t v_{p\lambda}}{v_{p\lambda}}_\pi=\ip{ v_{\lambda}^p}{v_{p\lambda}}_\pi<\infty$
and $g_t(y)\rightarrow \ip{ v_{\lambda}^p}{v_{p\lambda}}_\pi\ip{v_{p\lambda}}{v_{p\lambda}}_\pi^{-1}<\infty$ as $t\rightarrow\infty$ for all $y\in\Real$.

We also note that since $g_t\in L^1(\pi_{p\lambda})$, we have $f_t\in L^1(\tilde \pi_{p\lambda})$
where $\tilde\pi_{\mu}:=\ip{1}{ v_{\mu}}_\pi^{-1} v_\mu \pi $ and then
\[
\int_{y\in\Real} \tilde\pi_{p\lambda}(y) f_t(x,y)dy
= e^{p\lambda x}
\frac{\ip{v_\lambda^p}{ v_{p\lambda}}_\pi}{\ip{1}{ v_{p\lambda}}_\pi}\,
 e^{-(pE_\lambda - E_{p\lambda}) t}.
\]

\bigskip\noindent\textbf{The sum term.}
Note that under the parameter assumptions we have $\ip{ R v_{\lambda}^p}{v_{p\lambda}}_\pi<\infty$.
As for the finite-type model the fission times $S_u$ on the spine
occur as a Cox process and therefore
\begin{align*}
g_t(x,y) &:= e^{\lambda x}v_\lambda(y)
\tilde{\Q}_\lambda^{x,y} \bseq{\sum_{u<\xi_t}
v_\lambda(\eta_{S_u})^q e^{q\lambda \xi_{S_u}- qE_\lambda S_u}}\\
&= e^{\lambda x}v_\lambda(y)\int_0^t \tilde{\Q}_\lambda^{x,y} \bseq{2 R(\eta_s) \,
{v_\lambda(\eta_s)^q
e^{q\lambda \xi_s- qE_\lambda s}}} \, \diffd s\\
&=e^{\lambda x}v_\lambda(y)\int_0^t \tilde\Q_{p\lambda}^{0,y} \bseq{2R(\eta_s) \,
\frac{v_\lambda^p}{v_{p\lambda}}(\eta_s)}
e^{-(pE_\lambda -E_{p\lambda}) s} \, \diffd s \\
&= e^{p\lambda x} v_{p\lambda}(y)
k_t(y)
\end{align*}
where
\[
k_t(y):=\int_0^t
h_s(y)
 e^{-(pE_\lambda -E_{p\lambda}) s} \, \diffd s,
\qquad
h_s(y):=
\tilde\Q_{p\lambda}^{0,y} \bseq{2R(\eta_s) \,\frac{v_\lambda^p}{v_{p\lambda}}(\eta_s)}
\]
and $h_t,k_t\in L^1(\pi_{p\lambda})$. Note, $k_t(y)\uparrow k_\infty(y)\in L^1(\pi_{p\lambda})$
as $t\rightarrow\infty$ where
\[
\frac{\ip{k_t v_{p\lambda}}{v_{p\lambda}}_\pi}{\ip{v_{p\lambda}}{v_{p\lambda}}_\pi}
 = \ip{ 2 R v_{\lambda}^p}{v_{p\lambda}}_\pi
\frac{\left(1-e^{-(pE_\lambda -E_{p\lambda})t}\right)}{(pE_\lambda -E_{p\lambda})}
\uparrow
\frac{\ip{2 R v_{\lambda}^p}{v_{p\lambda}}_\pi} {(pE_\lambda -E_{p\lambda})}
= \ip{k_\infty v_{p\lambda}}{v_{p\lambda}}_\pi
<\infty.
\]
Note, $k_t\in L^1(\pi_{p\lambda})$ implies $g_t\in L^1(\tilde\pi_{p\lambda})$,
with an explicit calculation again possible.

\bigskip\noindent Bringing together the results for the sum and spine terms,
we have an upper bound
\begin{equation}\label{OU_L-p_bound}
P^{x,y}(Z_\lambda(t)^p) \leq e^{p\lambda x}
{v_{p\lambda}(y)}
\bbset{ g_t(y) e^{-(pE_\lambda - E_{p\lambda}) t} + k_t(y)}
\in L^1(\tilde \pi_{p\lambda})
\end{equation}
and hence the submartingale property reveals that
\[
P^{x,y}(Z_\lambda(t)^p) \leq e^{p\lambda x}
{v_{p\lambda}(y)}k_\infty(y)<\infty
\]
for all $t\geq0$ and all $y\in\Real$. \qed

\bigskip\noindent\textbf{Proof of Part 2:} As there is family sizes pose no problems,
we need only dominate the martingale by the spine at fixed time $t$:
\[
\tilde{\Q}_\lambda^{x,y}(Z_\lambda(t)^q)
\geq
\tilde{\Q}_\lambda^{x,y}(v_\lambda(\eta_t)^q e^{q\lambda\xi_t-qE_\lambda t})
=
e^{q\lambda x}\,\frac{v_{p\lambda}(y)}{v_\lambda(y)}
\, \tilde{\Q}^y_{p\lambda}\left(\frac{v_\lambda^p}{v_{p\lambda}}(\eta_t)\right)
e^{-(pE_\lambda -E_{p\lambda}) t}
\]
and $Z_\lambda$ is therefore unbounded in $\mathcal{L}^p(P^x)$ if either $pE_\lambda-E_{p\lambda}<0$
or $\ip{v_\lambda^p}{v_{p\lambda}}_\pi=\infty$.
\qed

\bigskip\noindent\textbf{Proof of Part 3:} The proof that we have
seen in the finite-type model will work here with little change:
under $\tilde{\Q}_\lambda$ the spatial motion is
\[
\xi_t =  B\left(\int_0^t a(\eta_s) \diffd s\right) + \lambda
\int_0^t a(\eta_s) \diffd s,
\]
and the type process $\eta_s$ has invariant distribution $N(0,\frac{\theta}{2\mu_\lambda})$, whence
$t^{-1}\xi_t \to \lambda a \theta/ \mu_\lambda = E_\lambda'$.
Therefore it follows that under $\tilde{\Q}_\lambda$ the diffusion
$\xi_t+c_\lambda t$ drifts off to $-\infty$ if
$\lambda\outsidewithout$.
When $\lambda=\tilde\lambda$, it is also simple to check that
$\xi_t+c_\lambda t$ is recurrent, so has $\liminf\{\xi_t+c_\lambda t\}=-\infty$.
Whence, in either case, bounding $Z_\lambda$ below by the spine's contribution as done before,
we have $Z_\lambda(t)\geq v_\lambda(\eta_t)e^{\lambda(\xi_t+c_\lambda t)}$
and since $v_\lambda>0$ and $\eta_t$ recurrent, we see that
$\limsup_{t\to\infty} Z_\lambda(t) =\infty$
almost surely under $\tilde{\Q}_\lambda^{x,y}$. \qed

\def\polhk#1{\setbox0=\hbox{#1}{\ooalign{\hidewidth
  \lower1.5ex\hbox{`}\hidewidth\crcr\unhbox0}}} \def\cprime{$'$}
\providecommand{\bysame}{\leavevmode\hbox
to3em{\hrulefill}\thinspace}
\providecommand{\MR}{\relax\ifhmode\unskip\space\fi MR }
\providecommand{\MRhref}[2]{%
  \href{http://www.ams.org/mathscinet-getitem?mr=#1}{#2}
} \providecommand{\href}[2]{#2}

\end{document}